\documentclass{article}

\usepackage{amsmath}
\usepackage{latexsym}
\usepackage{graphicx}
\usepackage{psfrag}
\usepackage{amssymb}
\usepackage{amscd}

\def\<{\langle}
\def\>{\rangle}
\def\0{{{\bf 0}}}
\def\calC{{\mathcal C}}
\def\CD{{\mathcal D}}
\def\OO{{\mathcal O}}
\def\CE{{\mathcal E}}
\def\CF{{\mathcal F}}
\def\tCD{{\tilde {\mathcal D}}}

\def\CS{{\mathbb S}}
\def\CT{{\mathcal T}}
\def\CA{{\mathcal A}}
\def\CB{{\mathcal B}}
\def\CM{{\mathcal M}}

\def\CV{{\mathcal V}}

\def\calC{{\mathcal C}}

\def\AA{{\mathbb A}}
\def\CC{{\mathbb C}}

\def\FF{{\mathbb F}}

\def\II{{\mathbb I}}

\def\PP{{\mathbb P}}
\def\QQ{{\mathbb Q}}
\def\RR{{\mathbb R}}
\def\TT{{\mathbb T}}
\def\ZZ{{\mathbb Z}}

\def\pp{{\mathfrak p}}

\def\tX{{\tilde X}}

\def\tA{{\tilde A}}
\def\tB{{\tilde B}}

\newcommand{\spec}{\operatorname{Spec}}

\newcommand{\im}{\operatorname{im}}

\newcommand{\Hom}{\operatorname{Hom}}

\newcommand{\rank}{\operatorname{rank}}
\newcommand{\End}{\operatorname{End}}
\newcommand{\Fr}{\operatorname{Fr}}
\newcommand{\Lie}{\operatorname{Lie}}
\newcommand{\Fabs}{F_{\mbox{\rm \tiny abs}}}

\newcommand{\gal}{\operatorname{Gal}}
\newcommand{\Aut}{\operatorname{Aut}}
\newcommand{\DR}{\operatorname{DR}}
\newcommand{\GL}{\operatorname{GL}}
\newcommand{\Sp}{\operatorname{Sp}}
\newcommand{\GU}{\operatorname{GU}}

\newcommand{\univ}{\mbox{\rm \tiny univ}}

\newcommand{\pnew}{\mbox{\rm \tiny $p$-new}}

\newcommand{\red}{\mbox{\rm \tiny red}}

\newcommand{\alg}{\mbox{\rm \tiny alg}}
\newcommand{\cris}{\mbox{\rm \tiny cris}}
\newcommand{\smallss}{\mbox{\rm \tiny ss}}
\newcommand{\et}{\mbox{\rm \tiny \'et}}

	{\begin{list}
		{\noindent\makebox[0mm][r]{\arabic{enumi}.}}
		{\usecounter{enumi} \topsep=1.5mm \itemsep=0mm}
	}
	{\end{list}}

	{\begin{list}
		{\noindent\makebox[0mm][r]{\arabic{enumi}.}}
		{\usecounter{enumi} \topsep=1.5mm \itemsep=-1mm}
	}
	{\end{list}}

\begin{document}

  \newtheorem{theorem}{Theorem}[section]
  \newtheorem{definition}[theorem]{Definition}
  \newtheorem{lemma}[theorem]{Lemma}
  \newtheorem{corollary}[theorem]{Corollary}
  \newtheorem{proposition}[theorem]{Proposition}
  \newtheorem{conjecture}[theorem]{Conjecture}
  \newtheorem{question}[theorem]{Question}
  \newtheorem{problem}[theorem]{Problem}

  \newtheorem{example}[theorem]{Example}
  \newtheorem{remark}[theorem]{Remark}

\newenvironment{proof}{{{\noindent{\it Proof.\ }}}}{{\hfill $\Box$}}

\title{A geometric Jacquet-Langlands correspondence for $U(2)$ Shimura
varieties}  

\author{David Helm
\footnote{University of Texas at Austin; dhelm@math.utexas.edu}}

\maketitle

Let $G$ be a unitary group over $\QQ$, associated to a CM-field $F$
with totally real part $F^+$, with signature $(1,1)$
at all the archimedean places of $F^+$.  Under certain hypotheses
on $F^+$, we show that Jacquet-Langlands correspondences
between certain automorphic representations of $G$ and
representations of a group $G^{\prime}$ isomorphic to $G$
except at infinity can be realized in the cohomology of Shimura
varieties attached to $G$ and $G^{\prime}$.

We obtain these Jacquet-Langlands correspondences by
studying the bad reduction of a Shimura variety $X$ attached
to $G$ at a prime $p$ for which $X$ has ``$\Gamma_0(p)$''
level structure, and construct a ``Deligne-Rapoport'' model
for $X$.  The irreducible components of the special fiber of
this model
have a global structure that can be explicitly described in
terms of Shimura varieties $X^{\prime}$ for unitary
groups $G^{\prime}$ isomorphic to $G$ away from infinity.

The weight spectral sequence of Rapoport-Zink then yields
an expression for certain pieces of the weight filtration
on the {\'e}tale cohomology of $X$ in terms of the cohomology
of $X^{\prime}$.  This identifies a piece of this weight
filtration with a space of algebraic modular forms for
$G^{\prime}$.  This result implies certain cases of the 
Jacquet-Langlands correspondence for $G$ and $G^{\prime}$ in terms of
a canonical map between spaces of arithmetic interest, rather than 
simply as an abstract bijection between isomorphism classes of 
representations.


\section{Introduction}

Let $G$ and $G^{\prime}$ be algebraic groups over $\QQ$, isomorphic over
all non-archimedean primes of $\QQ$ but not necessarily isomorphic over
$\RR$.  In many cases, the Jacquet-Langlands correspondence predicts
an bijection between automorphic representations of $G$ and $G^{\prime}$
satisfying certain technical conditions.  (See, for instance,~\cite{JL}
for the case of quaternionic algebraic groups.)

Traditionally, Jacquet-Langlands correspondences for various $G$ and 
$G^{\prime}$ are established via trace calculations in the theory
of automorphic representations.  One uses these calculations to
show that certain spaces of representations are abstractly isomorphic,
and deduces a Jacquet-Langlands correspondence from the existence
of such an isomorphism.

In the case of $\GL_2$, an observation of Serre \cite{Se}
shows an alternative way of proceeding:  Serre considers spaces of
modular forms as sections of line bundles on a modular curve.  By restricting
these line bundles to the supersingular locus of the characteristic
$p$ fiber of such a curve, Serre relates spaces of modular forms mod $p$
to spaces of functions on this supersingular locus.  This
supersingular locus is a double coset space for a quaternion
algebra $B$, and one can interpret the space of functions on such
a space as a space of ``mod $p$ modular forms'' for $B^{\times}$. 

In contrast to the traditional approach, which yields only a bijection
between isomorphism classes of representations, Serre's result gives a
canonical isomorphism between spaces that arise naturally from the
geometry of the Shimura varieties attached to the groups $\GL_2$
and $B^{\times}$.  

More recently, Ghitza \cite{Gh1},~\cite{Gh2} has generalized Serre's 
techniques to the case of Siegel modular varieties; he connects
automorphic forms for $\Sp_{2g}$ to algebraic modular forms~\cite{gross} 
on a quaternion unitary group.  As his results
(like those of Serre) use mod $p$ coefficients rather than complex
coefficients, his work has consequences for congruences of
Siegel modular forms; in particular he can show that all Siegel Hecke
eigensystems are cuspidal mod $p$. 

One disadvantage of this approach is that it only establishes a
correspondence in characteristic $p$.  An alternative approach, that
also establishes Jacquet-Langlands correspondences via the geometry
of Shimura varieties, is that of Ribet~\cite{ribet}.  Ribet considers
character groups of tori attached to the bad reduction of modular curves
and Shimura curves at various primes; these are free $\ZZ$-modules
with actions of Hecke operators for $\GL_2$ or a quaternion algebra $B$.
He then gives an exact sequence relating these character groups; this
amounts to a Jacquet-Langlands correspondence between $\GL_2$ and
$B^{\times}$ with {\em integer} coefficients.  As with the approach
of Serre and Ghizta, Ribet's construction is canonical and contains
data about congruences (in fact, it is a key ingredient of his proof
of level-lowering.)

Here we obtain similar results in the case of a Shimura variety $X$ attached
to a unitary group $G$ whose real points are (up to a factor of $\RR^{\times}$)
are isomorphic to a product of $U(1,1)$'s.  Such $G$ are defined in terms of a 
CM field $F$, with totally real part $F^+$.  Our approach most closely mirrors
that of Ribet, in that it makes use of semistable models of these varieties
at primes of bad reduction in order to deduce results that hold in characteristic
zero (as opposed to the mod $p$ results of Serre and Ghitza, which hold for
primes of good reduction.)  The reduction theory of Jacobians that plays a prominent
role in Ribet's theory is replaced with the theory of vanishing cycles and weight
filtrations for these more general Shimura varieties.

Making this argument work requires a through understanding of an analogue
of the ``Deligne-Rapoport model'' for the Shimura varieties under consideration.
The local structure of such models is by now very well-understood~\cite{RZ1}.
Unfortunately, cohomological arguments of the sort we have in mind seem to require an
understanding of the {\em global} structure of various strata in such a model; questions
of this nature have recieved far less attention.  

The bulk of the paper 
(sections~\ref{sec:stratify} through~\ref{sec:DR}) is devoted to establishing the
necessary global results.  In section~\ref{sec:stratify} we define an
``$a$-number stratification'' of the Shimura varieties under consideration at
a prime of good reduction, that is exactly analogous to the stratification of
Hilbert modular varieties introduced by E. Goren and F. Oort in~\cite{stratification},)  
and use standard techniques to understand the local structure of the strata.
In section~\ref{sec:sparse} we study the global structure of some of these strata;
surprisingly, many of the strata that arise in this way turn out to be (nearly)
$(\PP^1)^r$-bundles over Shimura varieties attached to other algebraic groups 
$G^{\prime}$.
These groups are isomorphic to $G$ at all finite places but not at infinity.  This
is already highly suggestive of a ``Jacquet-Langlands'' correspondence.  Section
\ref{sec:DR} gives a complete, {\em global} description of the components of
the special fiber of $X$ at certain
primes of bad reduction; the result gives each component of the special
fiber as (nearly) a $(\PP^1)^r$-bundle over one of the Shimura varieties
attached to the groups $G^{\prime}$ arising in the previous section.

Section~\ref{sec:cohomology} uses these calculations to deduce the main arithmetic
results of the paper.  The key tool is the Rapoport-Zink weight spectral
sequence~\cite{RZ2},~\cite{Saito}, which computes the weight filtration on
the {\'e}tale cohomology of an algebraic variety with semistable reduction
in terms of the cohomology of strata of the special fiber.  We apply this
to a semistable resolution of the model for $X$ described in
section~\ref{sec:DR}.
The result is Theorem~\ref{thm:filtration}, which gives an isomorphism
between the highest weight quotient of the
the middle {\'e}tale cohomology of $X$ and a space of
algebraic modular forms for a group $G^{\prime}$ that is isomorphic to $G$
at all finite places but (unlike $G$) is compact at infinity.  This isomorphism
is Hecke-equivariant, and therefore provides an explicit, canonical geometric
realization of Jacquet-Langlands transfers between $G$ and $G^{\prime}$.  For
surfaces, the isomorphism is even integral; that is, it works with $\ZZ_l$
coefficients as well as $\QQ_l$ coefficients.

There are some notational conventions we use throughout the paper.  First,
we often use the same symbol for an isogeny of abelian varieties, and for
the map on de Rham or crystalline cohomology it induces.  The dual of
an abelian variety $A$ will be denoted $A^{\vee}$.

\section{Basic definitions and properties} \label{sec:basic}

We begin with the definition and basic properties of $U(2)$ Shimura varieties.

Fix a totally real field $F^+$, of degree $d$ over $\QQ$.  Let $E$ be an
imaginary quadratic extension of $\QQ$, of discriminant $N$, and let $x$ be
a square root of $D$ in $E$.  Let $F$ be the field $EF^+$.

Fix a square root $\sqrt{D}$ of $D$ in $\CC$.  Then any embedding
$\tau: F^+ \rightarrow \RR$ induces two embeddings 
$p_{\tau}, q_{\tau}: F \rightarrow \CC$, via
\begin{eqnarray*}
p_{\tau}(a+bx) & = & \tau(a) + \tau(b)\sqrt{N} \\
q_{\tau}(a+bx) & = & \tau(a) - \tau(b)\sqrt{N}.
\end{eqnarray*}

Fix a two-dimensional $F$-vector space $\CV$, equipped with an alternating, 
nondegenerate
pairing $\langle, \rangle: \CV \times \CV \rightarrow \QQ.$  We require that
$$\langle \alpha x,y\rangle = \langle x, \overline{\alpha} y \rangle$$ 
for all $\alpha$ in $F$.  (Here $\overline{\alpha}$ denotes the complex
conjugate of $\alpha$.)

Each embedding $\tau: F^+ \rightarrow \RR$ gives us a complex vector
space $\CV_{\tau} = \CV \otimes_{F^+,\tau} \RR$.  The pairing $\langle,\rangle$
on $\CV$ induces a pairing $\langle,\rangle_{\tau}$ on $\CV_{\tau}$; 
this pairing is the ``imaginary part'' of a unique Hermitian pairing
$[,]_{\tau}$ on $\CV_{\tau}$.  We denote the
number of $1$'s in the signature of $[,]_{\tau}$ by $r_{\tau}(\CV)$, and the
number of $-1$'s by $s_{\tau}(\CV)$.  If $\CV$ is obvious from the context,
we will often omit it, and denote $r_{\tau}(\CV)$ and $s_{\tau}(\CV)$
by $r_{\tau}$ and $s_{\tau}$.

We fix a $\hat \ZZ$-lattice $\CT$ in $\CV(\AA^f_{\QQ})$, stable under
the action of $\OO_F$, such that the pairing $\langle,\rangle$
induces a $\hat \ZZ$-valued pairing on $\CT$.

Finally, we choose a prime $p \in \QQ$, unramified in $F$, and split in $E$, 
such that the cokernel
of the induced map $\CT \rightarrow \Hom(\CT,\hat \ZZ)$ has order prime to $p$.

Let $k_0$ be a finite field of characteristic $p$ large enough to contain
all the residue fields of $\OO_F$ of characteristic $p$.  We fix an
embedding of $W(k_0)$ into $\CC$.  Such a choice associates to
each $p_{\tau}$ (resp. $q_{\tau}$) as above a map 
$\OO_F \rightarrow W(k_0)$, and therefore also a map
$\OO_F/p\OO_F \rightarrow k_0$.  In an abuse of notation
we denote each of these by $p_{\tau}$ (resp. $q_{\tau}$), and rely on
the context to make it clear exactly which we mean.  

All these choices allow us to define the algebraic group associated to
the Shimura varieties of interest to us.  Specifically, we let
$G$ be the algebraic group over $\QQ$ such that for any $\QQ$-algebra
$R$, $G(R)$ is the subgroup of $\Aut_F(\CV \otimes_{\QQ} R)$ consisting
of all $g$ such that there exists an $r$ in $R^{\times}$ with
$\langle gx,gy \rangle = r\langle x,y \rangle$ for all $x$ and $y$ in
$\CV \otimes_{\QQ} R$.  In particular, $G(\RR)$ is the subgroup of
$$\prod_{\tau: F^+ \rightarrow \RR} \GU(r_{\tau},s_{\tau})$$
consisting of those tuples $(g_{\tau})_{\tau: F^+ \rightarrow \RR}$
such that the ``similitude ratio'' of $g_{\tau}$ is the same for all
$\tau$.

Now fix a compact open subgroup $U$ of $G(\AA^f_{\QQ})$, such that
$U$ preserves $\CT$, and such that $U_p$ is the 
largest subgroup of $G(\QQ_p)$ preserving $\CT_p$.  If $U$ is
sufficiently small, there exists a scheme $X_U(\CV)$ over $W(k_0)$
representing the functor that takes an $W(k_0)$-scheme $S$ to the set
of isomorphism classes of triples $(A,\lambda, \rho)$ where
\begin{enumerate}
\item $A$ is an abelian scheme over $S$ of dimension $2d$, with an action
of $\OO_F$, such that the characteristic polynomial of $\alpha \in \OO_F$ 
as an endomorphism of $\Lie(A/S)$ is given by
$$\prod_{\tau: F^+ \rightarrow \RR} (x - p_{\tau}(\alpha))^{r_{\tau}(\CV)} 
(x - q_{\tau}(\alpha))^{s_{\tau}(\CV)}.$$
\item $\lambda$ is a polarization of $A$, of degree prime to $p$, such
that the Rosati involution associated to $\lambda$ induces complex
conjugation on $\OO_F \subset \End(A)$.
\item $\rho$ is a $U$-orbit of isomorphisms 
$T_{\hat \ZZ^{(p)}} A \rightarrow \CT^{(p)}$, sending the 
Weil pairing on 
$T_{\hat \ZZ^{(p)}} A$
to a scalar multiple of the pairing $\langle, \rangle$ on 
$\CT^{(p)}$.  (Here
$T_{\hat \ZZ^{(p)}} A$ denotes the product over all $l \neq p$ of the $l$-adic 
Tate modules of $A$.)
\end{enumerate}
We call a scheme $X_U(\CV)$ arising in this way a $U(2)$ Shimura variety.

\begin{remark} \rm \label{rem:decomp}
The condition on the characteristic polynomials of
elements of $\OO_F$ acting on $\Lie(A/S)$ can be reformulated in
a more illuminating manner.  If $M$ is a locally free sheaf on
$S$ with an action of $\OO_F$, then $M$ admits an action of
$\OO_F \otimes W(k_0)$, and hence breaks up as a direct sum:
$$M = \bigoplus_{\varphi: \OO_F \rightarrow W(k_0)} M_{\varphi},$$
where $M_{\varphi}$ is the locally free submodule of $M$ on
which $\OO_F$ acts via the embedding $\varphi$ of $\OO_F$ in $W(k_0)$.  Our choice
of embedding of $W(k_0)$ in $\CC$ allows us to consider each
complex embedding $p_{\tau}$ or $q_{\tau}$ of $\OO_F$ as such
a character.  In this language, the characteristic polynomial condition is
equivalent to requiring that $\Lie(A/S)_{p_{\tau}}$ and $\Lie(A/S)_{q_{\tau}}$
are locally free of ranks $r_{\tau}(\CV)$ and $s_{\tau}(\CV)$,
respectively.
\end{remark}

Since we have taken $U_p$ to be the subgroup of $G(\QQ_p)$ preserving
$\CT_p$, $U_p$ is a hyperspecial subgroup of $G(\QQ_p)$.  It follows
that the mod $p$ reduction $X_U(\CV)_{/k_0}$ is a smooth variety.
Its dimension is equal to the number of $\tau$ such that 
$r_{\tau}(\CV) = s_{\tau}(\CV) = 1$, or equivalently to the sum
$$\sum_{\tau} r_{\tau}(\CV)s_{\tau}(\CV).$$

We will also be interested in a variant of $X_U(\CV)$ that has some
level structure at $p$.  In particular, we define $X_{U,p}(\CV)$
to be the $W(k_0)$-scheme representing the functor that takes 
$S$ to the set of isomorphism classes of tuples 
$(A, \lambda, \rho, B, \lambda^{\prime},\phi)$,
where
\begin{enumerate} 
\item $(A, \lambda, \rho)$ is an $S$-point of $X_U(\CV)$.
\item $B$ is an abelian scheme over $S$ of dimension $2d$, with an action
of $\OO_F$.
\item $\lambda^{\prime}$ is a polarization of $B$, of degree prime to $p$.
\item $\phi: A \rightarrow B$ is an $\OO_F$-isogeny of degree $p^{2d}$, 
satisfying:
\begin{itemize}
\item $p\lambda = \phi^{\vee}\lambda^{\prime}\phi$,
\item for each prime $\pp$ of $\OO_F$ over $p$, $\ker \phi \cap A[\pp]$
is a finite flat group scheme over $S$ of order equal to the norm of $\pp$, and
\item the cokernel of the map $H^1_{\DR}(B/S) \rightarrow H^1_{\DR}(A/S)$
is a locally free $\OO_F \otimes \OO_S/p$-module of rank one.
\end{itemize}
\end{enumerate}

Unlike $X_U(\CV)$, $X_{U,p}(\CV)$ is not smooth; in fact,
in characteristic $p$ it is not even irreducible.  The structure of its 
irreducible components will be discussed in section~\ref{sec:DR}.

\begin{remark} \rm The schemes $X_U(\CV)$ and $X_{U,p}(\CV)$
depend not only on $U$ and $\CV$ but on all of the choices we have
made in this section.  To avoid clutter, we have chosen to suppress 
most of these choices in our notation.
\end{remark}
 
One can also consider a more general version of the above moduli problems.
Let $F^+$, $E$, and $F$ be as above, and let $D$ be a central simple algebra 
over $F$ of dimension $4$.  Fix an involution 
$\alpha \mapsto \overline{\alpha}$ of $D$, of the second kind, and a 
maximal order $\OO_D$ of $D$ stable under this involution.

Let $\CV$ now denote a free left $D$-module of rank one, equipped with
a nondegenerate alternating pairing 
$$\langle,\rangle: \CV \times \CV \rightarrow \QQ$$
with the property that
$$\langle \alpha x,y \rangle = \langle x, \overline{\alpha} y\rangle$$
for all $\alpha$ in $D$.  For each $\tau: F^+ \rightarrow \RR$,
$\CV \otimes_{F^+,\tau} \RR$ is a free $D \otimes_{F^+,\tau} \RR$-module
of rank one.  Since $F$ is purely imaginary, $D \otimes_{F^+,\tau} \RR$
is isomorphic to $\GL_2(\CC)$, so choosing a suitable idempotent of
$e$ of $D \otimes_{F^+,\tau} \RR$, such that $e = \overline{e}$,
we find that $e(\CV \otimes_{F^+,\tau} \RR)$ is a two-dimensional
complex vector space with equipped with a nondegenerate Hermitian pairing
induced by $\langle,\rangle$.  We let $r_{\tau}(\CV)$ and $s_{\tau}(\CV)$
be the number of $1$'s and $-1$'s, respectively, in the signature of this 
pairing.

Finally, we fix a $\hat \ZZ$-lattice $\CT$ in $\CV(\AA_{\QQ}^f)$, stable under 
the action of $\OO$, such that $\langle, \rangle$ induces a $\hat \ZZ$-valued 
pairing on $\CT$.  We fix a $p$, unramified in $F^+$ and split in $E$, such
that $B$ is split at $p$ and the cokernel of the map 
$\CT \rightarrow \Hom(\CT,\hat \ZZ)$ induced by $\langle,\rangle$ has order
prime to $p$.

Then for $G \subset \Aut_D(\CV)$ the subgroup of automorphisms which
take $\langle,\rangle$ to a scalar multiple of itself, and $U$
a compact open subgroup of $G$, preserving $\CT$ such 
that $U_p$ is the subgroup of
$G_p$ fixing $\CT_p$, we can associate a Shimura variety $X^D_U(\CV)$.
The corresponding moduli problem over $W(k_0)$ parametrizes tuples
$(A,\lambda,\rho)$, where
\begin{enumerate}
\item $A$ is an abelian scheme over $S$ of dimension $4d$, with an action
of $\OO_D$, such that the characteristic polynomial of $\alpha \in \OO_F$ 
as an endomorphism of $\Lie(A/S)$ is given by
$$\prod_{\tau: F^+ \rightarrow \RR} (x - p_{\tau}(\alpha))^{2r_{\tau}(\CV)} 
(x - q_{\tau}(\alpha))^{2s_{\tau}(\CV)}.$$
\item $\lambda$ is a polarization of $A$, of degree prime to $p$, such
that the Rosati involution associated to $\lambda$ induces the involution
$\alpha \mapsto \overline{\alpha}$ on $\OO_D$. 
\item $\rho$ is a $U$-orbit of isomorphisms 
$T_{\hat \ZZ^{(p)}} A \rightarrow \CT^{(p)}$, sending the 
Weil pairing on 
$T_{\hat \ZZ^{(p)}} A$
to a scalar multiple of the pairing $\langle, \rangle$ on 
$\CT^{(p)}$.  (Here
$T_{\hat \ZZ^{(p)}} A$ denotes the product over all $l \neq p$ of the $l$-adic 
Tate modules of $A$.)
\end{enumerate}

Similarly, one defines $X^D_{U,p}(\CV)$ to be the space that parameterizes
tuples $(A,\lambda,\rho,B,\lambda^{\prime},\phi)$, where
\begin{enumerate} 
\item $(A, \lambda, \rho)$ is an $S$-point of $X^D_U(\CV)$.
\item $B$ is an abelian scheme over $S$ of dimension $4d$, with an action
of $\OO_D$.
\item $\lambda^{\prime}$ is a polarization of $B$, of degree prime to $p$.
\item $\phi: A \rightarrow B$ is an $\OO_F$-isogeny of degree $p^{4d}$, 
satisfying:
\begin{itemize}
\item $p\lambda = \phi^{\vee}\lambda^{\prime}\phi$,
\item for each prime $\pp$ of $\OO_F$ over $p$, $\ker \phi \cap A[\pp]$
is a finite flat group scheme over $S$ of order equal to the norm of $\pp$ 
squared, and
\item the cokernel of the map $H^1_{\DR}(B/S) \rightarrow H^1_{\DR}(A/S)$
is a locally free $\OO_F \otimes \OO_S/p$-module of rank two.
\end{itemize}
\end{enumerate}

Note that for $D = M_2(F)$, with involution equal to conjugate transpose,
this is equivalent to our original moduli problem.
Specifically, given a point $(A,\lambda,\rho)$ of $X_U(\CV)$, 
we obtain a point $(A^2,\lambda^2,\rho^2)$ of $X_U^D(\CV)$ by letting
$\GL_2(\OO_F)$ act on $A^2$ in the obvious way.  Conversely, given
a point $(A,\lambda,\rho)$ of $X_U^D(\CV)$, for a suitable
idempotent $e$ of $M_2(\OO_F)$ (i.e., one stable under conjugate transpose)
$\lambda$ induces a prime-to-$p$ polarization 
$e\lambda: eA \rightarrow (eA)^{\vee}$,
and so $(eA,e\lambda,e\rho)$ is a point of $X_U(\CV)$. 

For the purposes of most of this paper, it will matter little whether
we work with $X_U(\CV)$ or $X_U^D(\CV)$ for some $D$.  In particular, all
of our constructions are done in terms of the $p$-divisible groups
of the moduli objects we consider.  Since $D$ is split at $p$,
we can fix an idempotent $e$ of $\OO_D \otimes \ZZ_p$ such that
$e = \overline{e}$.  Then if we have a point $(A,\lambda,\rho)$ of 
$X_U^D(\CV)$, $eA[p^{\infty}]$ behaves exactly like the $p$-divisible
group of a moduli object of $X_U(\CV)$.

Therefore, to avoid cluttering our notation with idempotents everywhere,
for the bulk of the paper we consider only the schemes $X_U(\CV)$ and
$X_{U,p}(\CV)$.  All of our results hold more generally
for any $X_U^D(\CV)$ for $D$ as above (in particular, split at $p$).
We leave it to the reader to determine how the proofs given here may be
adapted to this case. 

\section{A stratification of $U(2)$ Shimura varieties mod $p$} 
\label{sec:stratify}

Fix choices
for $F^+$, $E$, $\CV$, and $U$, and let $r_{\tau} = r_{\tau}(\CV)$ and
$s_{\tau} = s_{\tau}(\CV)$ for all $\tau$.  Let 
$X = X_U(\CV)_{k_0}$ be
the fiber over $k_0$ of the 
$U(2)$ Shimura variety arising from these choices as in section~\ref{sec:basic}.
Our first goal will be to construct a stratification on $X$.  We do
this by adapting a closely related construction of Goren and 
Oort~\cite{stratification} for Hilbert modular varieties.  The
stratification we obtain is a ``generalized Ekedahl-Oort'' stratification
for $U(2)$ Shimura varieties; these have been studied in considerable
generality by many authors; see for instance~\cite{oort}, ~\cite{moonen},
or~\cite{wedhorn}.

Let $\CA$ be the universal abelian variety over $X$.  Then the relative
Lie algebra $\Lie(\CA/X)$ is a locally free $\OO_X$-module of rank $2d$,
and decomposes as a direct sum:
$$\Lie(\CA/X) = \bigoplus_{\tau \in \CS }
(\Lie(\CA/X)_{p_{\tau}} \oplus \Lie(\CA/X)_{q_{\tau}}),$$
where $\CS = \Hom_{\alg}(F^+,\RR)$
as in remark~\ref{rem:decomp}.  Moreover, the bundles $\Lie(\CA/X)_{p_{\tau}}$ and 
$\Lie(\CA/X)_{q_{\tau}}$ are locally free of rank $r_{\tau}$ and $s_{\tau}$, 
respectively.

Fix a $\tau$ such that $r_{\tau} = s_{\tau} = 1$.
The Verschiebung $V: \CA^{(p)} \rightarrow \CA$ induces maps
$$\Lie(\CA^{(p)}/X)_{p_{\tau}} \rightarrow \Lie(\CA/X)_{p_{\tau}}$$
for all $\tau$.  Moreover, there is a natural isomorphism:
$$\Lie(\CA^{(p)}/X) \cong \Fabs^*(\Lie(\CA/X)),$$ 
where $\Fabs$ is the absolute Frobenius on $X$.  
Thus we have a natural isomorphism:
$$\Lie(\CA^{(p)}/X)_{p_{\tau}} \cong (\Fabs^*(\Lie(\CA/X)))_{p_{\tau}}.$$

Let $\sigma: k_0 \rightarrow k_0$ be the
Frobenius automorphism of $k_0$.  Via our choice of embedding of $W(k_0)$
in $\CC$ in section~\ref{sec:basic}, an embedding
$\tau: F^+ \rightarrow \RR$ corresponds to a unique map
$\OO_{F^+} \rightarrow k_0$.  Composing this map with $\sigma$,
we obtain a new map $\OO_{F^+} \rightarrow k_0$.  
We let $\sigma\tau$ denote the map $F^+ \rightarrow \RR$ corresponding
to this new map.  Since the absolute
Frobenius induces $\sigma$ on $k_0$, we have
$$(\Fabs^*(\Lie(\CA/X)))_{p_{\tau}} \cong \Fabs^*(\Lie(\CA/X)_{p_{\sigma\tau}})$$
with this convention.

Suppose further that $r_{\sigma\tau} = s_{\sigma\tau} = 1$,
so that $\Lie(\CA/X)_{p_{\sigma\tau}}$ is a line bundle.
Pullback of a line bundle by Frobenius sends transition functions to
their $p$th powers, and therefore sends a line bundle to its $p$th tensor
power.  We thus obtain a canonical isomorphism
$$\Lie(\CA^{(p)}/X)_{p_{\tau}} \cong \Lie(\CA/X)_{p_{\sigma\tau}}^{\otimes p}.$$
In particular, for each $\tau$, Verschiebung yields a section $h_{\tau}$
of the line bundle $\Hom(\Lie(\CA/X)_{p_{\sigma\tau}}^{\otimes p},
\Lie(\CA/X)_{p_{\tau}})$.

\begin{remark} \rm These $h_{\tau}$ are close analogues of the ``partial Hasse
invariants'' constructed by Goren in the setting of Hilbert modular 
varieties~\cite{hasse}.
\end{remark}

\begin{remark} \rm We could just as easily have done this construction with the
$q_{\tau}$ instead of with the $p_{\tau}$.  The existence of a prime-to-$p$
polarization on $\CA$ means that the $p_{\tau}$ and $q_{\tau}$ components
of $\Lie(\CA/X)$ are dual to each other, and hence that the bundle and
section obtained from considering the $q_{\tau}$ components are canonically
isomorphic to those obtained from the $p_{\tau}$ components.
\end{remark}

\begin{proposition} Let $k$ be a perfect field of characteristic $p$,
$x \in X(k)$, and $A/k$ the abelian variety corresponding
to $x$.  Then $x$ lies in $Z(h_{\tau})$ if and only if the
$k$-vector space $\Hom(\alpha_p, A[p])_{p_{\tau}}$ has dimension at least
one.
\end{proposition}
\begin{proof}
Let $\CD$ denote the contravariant Dieudonn{\'e} module of $A[p]$.  The
action of $\OO_F$ on $\CD$ decomposes $\CD$ into a direct sum:
$$\CD = \bigoplus_{\tau} \CD_{p_{\tau}} \oplus \CD_{q_{\tau}}.$$
Let $F_{\CD}$ and $V_{\CD}$ denote the Frobenius and Verschiebung 
endomorphisms of $\CD$.  Then $\Hom(\alpha_p, A[p])_{p_{\tau}}$
is nonzero if and only if the quotient
$$(\CD/(\im F_{\CD} + \im V_{\CD}))_{p_{\tau}}$$ 
is nonzero.  This latter module can also be written as 
$$\CD_{p_{\tau}}/(F_{\CD}(\CD_{p_{\sigma^{-1}\tau}}) + 
V_{\CD}(\CD_{p_{\sigma\tau}})).$$

On the other hand, one has a canonical isomorphism
$$\CD \cong H^1_{\DR}(A/k),$$ 
that sends $V_{\CD}(\CD)$ to the subspace $\Lie(A/k)^*$ of $H^1_{\DR}(A/k)$.
In particular, $V_{\CD}(\CD_{p_{\tau}})$ has dimension 
$r_{\sigma^{-1}\tau} = 1$. Since the image of $V_{\CD}$
is equal to the kernel of $F_{\CD}$, we also have that 
$F_{\CD}(\CD_{p_{\sigma^{-1}\tau}})$
is one dimensional.  It follows that 
$\Hom(\alpha_p, A[p])_{p_{\tau}}$ is nonempty if and only if 
$V_{\CD}(\CD_{p_{\sigma\tau}})$ equals $F_{\CD}(\CD_{p_{\sigma^{-1}\tau}})$.

In terms of de Rham cohomology, this means that 
$\Hom(\alpha_p, A[p])_{p_{\tau}}$ is nonzero if and only if 
$\Lie(A/k)_{p_{\tau}}^*$ is contained in $\Fr(H^1_{\DR}(A^{(p)}/k))$,
where $\Fr$ is the relative Frobenius $A \rightarrow A^{(p)}$.  The latter 
is the kernel of the map induced by Verschiebung on de Rham cohomology.
Thus $\Hom(\alpha_p, A[p])_{p_{\tau}}$ is nonzero if and only if the
map 
$$h_{\tau}^*:\Lie(A)^*_{p_{\tau}} \rightarrow \Lie(A^{(p)})^*_{p_{\tau}}$$
induced by Verschiebung is the zero map, as required.  
\end{proof} 

For $S$ a subset of those $\tau$ with $r_{\tau} = r_{\sigma^{-1}\tau} = 1$, 
let $X_S$ denote the (scheme-theoretic) intersection of the $Z(h_{\tau})$ 
for $\tau$ in $S$.

\begin{proposition} \label{prop:smooth} 
If $X_S \subset X$ is nonempty, then $X_S$ is smooth, of codimension equal 
to the cardinality of $S$.  Moreover, if $S$ is nonempty, then $X_S$ is proper.
\end{proposition}

Before proving this we recall some basic facts about the deformation theory
of abelian varieties:

Let $S$ be a scheme, and $S^{\prime}$ a nilpotent
thickening of $S$ equipped with divided powers.  Let $\calC_{S^{\prime}}$
denote the category of abelian varieties over $S^{\prime}$, and
$\calC_S$ denote the category of abelian varieties over $S$.  For $A$ an
object of $\calC_{S^{\prime}}$, let $\overline{A}$ denote its base change  
to $\calC_S$.

Fix an $A$ in $\calC_{S^{\prime}}$, and consider the module 
$H^1_{\cris}(\overline{A}/S)_{S^{\prime}}$.  This is a locally free
$\OO_{S^{\prime}}$-module, and we have a canonical isomorphism:
$$H^1_{\cris}(\overline{A}/S)_{S^{\prime}} \cong H^1_{\DR}(A/S^{\prime}).$$
Moreover, we have a natural subbundle 
$$\Lie(A/S^{\prime})^* \subset H^1_{\DR}(A/S^{\prime}).$$
This gives us a subbundle of 
$H^1_{\cris}(\overline{A}/S)_{S^{\prime}}$ that lifts the subbundle 
$\Lie(\overline{A}/S)^*$ of $H^1_{\DR}(\overline{A}/S)$.

In fact, knowing this lift allows us to recover $A$ from $\overline{A}$.
More precisely, let $\calC_S^+$ denote the category of pairs 
$(\overline{A},\omega)$, where $\overline{A}$ is an object of $\calC_S$ and 
$\omega$ is a subbundle
of $H^1_{\cris}(\overline{A}/S)_{S^{\prime}}$ that lifts 
$\Lie(\overline{A}/S)^*.$  Then the construction outlined above gives
us a functor from $\calC_{S^{\prime}}$ to $\calC_S^+$.

\begin{theorem}[Grothendieck] \label{thm:grot}
The functor $\calC_{S^{\prime}} \rightarrow \calC_S^+$ 
defined above is an equivalence of categories.
\end{theorem}

\begin{proof} For an outline of the basic ideas, see~\cite{montreal}, pp.
116-118.  A complete proof can be found in~\cite{mazurmessing}, primarily
chapter II, section 1.
\end{proof}

We are now in a position to prove Proposition~\ref{prop:smooth}.  First
note that $X_S$ is defined by the vanishing of sections of $n$
line bundles, where $n$ is the cardinality of $S$.  Thus $X_S$ has
codimension at most $n$ in $X$.  It suffices to show
that the tangent space to $X_S$ at any $k$-point
$x = (A,\lambda,\rho)$ of $X_S$ ($k$ finite) 
has codimension $n$ in the tangent space to $X$.  

Let $\II$ be the scheme $\spec k[\epsilon]/\epsilon^2$.
The first-order deformations of $A$ are in bijection with
the lifts $\omega$ of $\Lie(A/k)^*$ to a subbundle of
$H^1_{\cris}(A/k)_{\II}$.  The $\OO_F$-action on
$A$ lifts to the first order deformation given by $\omega$ if and only if
$\omega$ is of the form
$$\omega = \bigoplus_{\tau} \omega_{p_{\tau}} \oplus \omega_{q_{\tau}},$$
where $\omega_{p_{\tau}}$ (resp. $\omega_{q_{\tau}}$) is a lift of
$\Lie(A/k)_{p_{\tau}}^*$ 
(resp. $\Lie(A/k)_{q_{\tau}}^*$) to
$(H^1_{\cris}(A/k)_{\II})_{p_{\tau}}$ (resp.
$(H^1_{\cris}(A/k)_{\II})_{q_{\tau}}.)$ 

The (prime-to-$p$) polarization $\lambda$ on $A$ induces a perfect pairing on
$H^1_{\cris}(A/k)_{\II}$.  Since the Rosati involution
induces complex conjugation on $\OO_F$, this restricts to a perfect
pairing:
$$ (H^1_{\cris}(A/k)_{\II})_{p_{\tau}} \times 
(H^1_{\cris}(A/k)_{\II})_{q_{\tau}} \rightarrow 
k[\epsilon]/\epsilon^2.$$ 
By Theorem~\ref{thm:grot}, $\lambda$ lifts to the first-order deformation
given by $\omega$ if and only if 
$\omega_{q_{\tau}} = \omega_{p_{\tau}}^{\perp}$ for all $\tau$, where
$\perp$ denotes orthogonal complement with respect to the pairing above.
(c.f.~\cite{sgee}, Corollary 2.2).

Finally, for any choice of $\omega$ for which the $\OO_F$-action and
$\lambda$ lift to the corresponding deformation, the level structure
$\rho$ lifts uniquely.  Thus specifying a tangent vector to $X$ at $x$ 
is equivalent to specifying, for each $\tau$, a lift of the
the $r_{\tau}$-dimensional subspace $\Lie(A/k)_{p_{\tau}}$ of 
the two dimensional vector space $H^1_{\DR}(A/k)_{p_{\tau}}$ to a
subbundle $\omega_{p_{\tau}}$ of
$(H^1_{\cris}(A/k)_{\II})_{p_{\tau}}.$  If $r_{\tau}$ is zero or
two, the lift is obviously unique.  For each $\tau$ with $r_{\tau} = 1$,
such lifts form a one-dimensional $k$-vector space.  

Let $v$ be a tangent vector to $X$ at $x$, corresponding to a lift
$(\hat A, \hat \lambda, \hat \rho)$ of $(A,\lambda,\rho)$.  For $\tau$
such that $r_{\tau} = r_{\sigma^{-1}\tau} = 1$, the vector
$v$ lies in $X_{\{\tau\}}$ if and only if the map
$$V: \Lie(\hat A/\II)^*_{p_{\tau}} \rightarrow 
\Lie(\hat A^{(p)}/\II)^*_{p_{\tau}}$$
is the zero map.  This in turn is true if and only if the lift
$\omega_{p_{\tau}}$ coming from $v$ is contained in the kernel of
the map
$$V: (H^1_{\cris}(A/k)_{\II})_{p_{\tau}} \rightarrow
(H^1_{\cris}(A^{(p)}/k)_{\II})_{p_{\tau}}.$$
Since both $\omega_{p_{\tau}}$ and this kernel are subbundles
of rank one, this condition determines $\omega_{p_{\tau}}$ uniquely. 

Thus to specify a tangent vector to $X_S$ at $x$, it suffices to
specify $\omega_{p_{\tau}}$ for $\tau$ {\em not} in $S$.  In particular
the tangent space to $X_S$ at $x$ has codimension $n$ in the tangent
space to $X$ at $x$, as required.

It remains to show that $X_S$ is proper for $S$ nonempty.
Let $R$ be a discrete valuation ring containing $k_0$,
let $K$ be its field of fractions, and let $k$ be its residue field.
Let $x: \spec K \rightarrow X_S$ be a point of $X_S$, and
let $A/K$ be the corresponding abelian variety.  Then there exists a finite
extension $K^{\prime}$ of $K$ such that $A_{K^{\prime}}$ has either good or
multiplicative reduction. 

Suppose it had multiplicative reduction, and let $\TT$ be the maximal
subtorus of the reduction.  The character group of $\TT$ is a free
$\ZZ$-module of rank at most $2d$ on which $\OO_F$ acts; this is only 
possible if $\TT$ had rank exactly $2d$, in which case $A_{K^{\prime}}$
has purely toric reduction.  But then Verschiebung would induce an
isomorphism on the Lie algebra of the special fiber, which is impossible
because it has a kernel on the general fiber.

Thus $A_{K^{\prime}}$ has good reduction.  In particular $A_{K^{\prime}}$ 
extends
to an abelian scheme over $\spec R^{\prime}$, where $R^{\prime}$ is
the integral closure of $R$ in $K^{\prime}$.  The polarization and level
structure on $A_{K^{\prime}}$ extend to $R^{\prime}$ as well.  We thus
obtain a map $\spec R^{\prime} \rightarrow X_S$ extending the map
$\spec K^{\prime} \rightarrow X_S$, so $X_S$ is proper.

\begin{remark} \rm It follows immediately from this result that
the nonempty $X_S$ are the closed strata in a stratification (the $a$-number
stratification) of $X$.
\end{remark}

\begin{lemma} If $r_{\tau} = s_{\tau} = 1$ for all $\tau$, then $X_S$ is 
nonempty for all $S$.
\end{lemma}
\begin{proof} It suffices to prove this for $S = \CS$.
Let $\CE$ be a supersingular
elliptic curve over $\overline{\FF}_p$, and let $A$ be the abelian variety
$\CE \otimes_{\ZZ} \OO_F$; that is, the abelian variety such that
$A(S) = \CE(S) \otimes_{\ZZ} \OO_F$ for all schemes $S$. 
Then both $T_{\hat \ZZ^{(p)}} A$ and
$\CT^{(p)}$ are free $\hat \OO_F^{(p)}$-modules of rank two.  
Fix any isomorphism
$\rho$ between them, compatible with the action of $\OO_F$.  The natural
polarization on $E$ then induces a polarization on $A$, that identifies
$T_{\hat \ZZ^{(p)}} A^{\vee}$ with $\Hom(\CT^{(p)}, \hat \ZZ^{(p)})$.
Thus the pairing on $\CT^{(p)}$ determines a
polarization $\lambda$ of $A$, of degree prime to $p$.  Moreover,
we have isomorphisms
$$\Lie(A/\overline{\FF}_p) \cong \Lie(\CE/\overline{\FF}_p) \otimes_{\ZZ}
\OO_F \cong \overline{\FF}_p \otimes_{\ZZ} \OO_F,$$
so the characteristic polynomial condition on $\Lie(A/\overline{\FF}_p)$ is 
satisfied.  The triple $(A,\lambda,\rho)$ thus determines a point of $X$.

We have isomorphisms
$$\Hom(\alpha_p,A[p]) \cong \Hom(\alpha_p,\CE[p]) \otimes_{\ZZ} \OO_F
\cong \overline{\FF}_p \otimes_{\ZZ} \OO_F,$$
and the $p_{\tau}$ component of the latter is nonzero for all $\tau$. 
Thus $\Hom(\alpha_p,A[p])_{p_{\tau}}$ is nonzero for all $\tau$,
so $(A,\lambda,\rho)$ lies on $X_{\CS}$.
\end{proof}

\section{Description of the sparse strata} \label{sec:sparse}

Now choose $\CV$ so that $r_{\tau}(\CV) = s_{\tau}(\CV) = 1$ for all 
$\tau$, and let $X = X_U(\CV)_{k_0}$.  
As we have seen, $X$ has a stratification 
whose closed strata $X_S$ are indexed by subsets $S$ of $\CS$.
We compute the global structure of $X_S$ for certain $S$ 
in terms of $U(2)$ Shimura varieties with different archimedean 
invariants.

\begin{definition} A subset $S$ of $\CS$
is {\em sparse} if for every $\tau$ in $S$, $\sigma\tau$ is not in $S$,
where $\sigma$ is the Frobenius automorphism of $k_0$.
\end{definition}

\begin{remark} \rm If $p$ splits completely in $F^+$, then there are no
sparse sets $S$, and the theory of this section is vacuous.  On
the other hand, if the residue class degree of
every prime of $F^+$ above $p$ is even (so that in particular $d$ is even),
then there exist sparse subsets of order up to $\frac {d}{2}.$
\end{remark}

Fix a sparse subset $S$.  Let $Y_S$ be the moduli space over $k_0$ 
whose points over $T$ are isomorphism
classes of
tuples $(A,\lambda,\rho,B,\lambda^{\prime},\phi)$, where
\begin{enumerate} 
\item $(A,\lambda,\rho)$ is a point of $X_U(T)$,
\item $B$ is an abelian scheme over $T$ of dimension $2d$, with an
action of $\OO_F$,
\item $\lambda^{\prime}$ is a polarization of $B$, of degree prime to $p$,
and
\item $\phi: A \rightarrow B$ is an $\OO_F$-isogeny of degree $p^2d$,
satisfying:
\begin{itemize}
\item $p\lambda = \phi^{\vee}\lambda^{\prime}\phi$,
\item for each $\tau$ not in $S$, the subbundle 
$\phi(H^1_{\DR}(B/T))_{p_{\tau}}$
of $H^1_{\DR}(A/T)_{p_{\tau}}$ is equal to $H^1_{\DR}(A/T)_{p_{\tau}}$,
and
\item for each $\tau$ in $S$, the subbundle $\phi(H^1_{\DR}(B/T))_{p_{\tau}}$
of $H^1_{\DR}(A/T)_{p_{\tau}}$
is {\em simultaneously} equal to both $\Lie(A/T)_{p_{\tau}}^*$ and
$\Fr(H^1_{\DR}(A^{(p)}/T)_{p_{\tau}})$, when considered as
subbundles of $H^1_{\DR}(A/T)_{p_{\tau}}$.  (Note that this implies that
$(A,\lambda,\rho)$ is actually a point of $X_S$.)
\end{itemize}
\end{enumerate}

The moduli space $Y_S$ can be thought of as parametrizing a point
$(A,\lambda,\rho)$ of $X_S$, plus an isogeny $\phi: A \rightarrow B$
that serves as a ``witness'' to the fact that $(A,\lambda,\rho)$
lies in $X_S$.  In fact, we will soon prove that the map
$Y_S \rightarrow X_S$ that forgets $\phi$ is an isomorphism.
First we need the following calculation:

\begin{lemma} \label{lemma:lierank}
The rank of $\Lie(B/T)_{p_{\tau}}$ is two if $\tau$ is in $S$,
zero if $\sigma\tau$ is in $S$, and one otherwise.  In all cases,
the rank of $\Lie(B/T)_{q_{\tau}}$ is equal to 
$2 - \rank \Lie(B/T)_{p_{\tau}}$.
\end{lemma}

\begin{proof}
Since $\Lie(B/T)_{p_{\tau}}$ is locally free, we can compute its rank
after base change to a geometric point of $T$.  We may thus assume
$T = \spec k$ for a perfect field $k$ containing $k_0$.
This allows us to apply Dieudonn{\'e} theory.

Let $\tCD_A$ (resp. $\tCD_B$) denote the Dieudonn{\'e} module of the 
$p$-divisible group $A[p^{\infty}]$ (resp. $B[p^{\infty}]$).  
Let $\CD_A$ (resp. $\CD_B$)
denote the Dieudonn{\'e} module of $A[p]$ (resp. $B[p]$), so that
$\CD_A \cong \tCD_A/p\tCD_A$ and similarly for $\CD_B$.

The space $\Lie(B/T)^*_{p_{\tau}}$ is isomorphic to the image of
$(\CD_B)_{p_{\sigma\tau}}$ in $(\CD_B)_{p_{\tau}}$ under Verschiebung.
Thus we have: 
\begin{eqnarray*}
\dim_k \Lie(B/T)_{p_{\tau}} & = & 
2 - \dim_k (\CD_B)_{p_{\tau}}/V_B((\CD_B)_{p_{\sigma\tau}}) \\
& = & 2 - \dim_k (\tCD_B)_{p_{\tau}}/V_B((\tCD_B)_{p_{\sigma\tau}})
\end{eqnarray*}
The map $\phi$ injects $\tCD_B$ into $\tCD_A$, so the latter is equal
to
$$2 - \dim_k \phi(\tCD_B)_{p_{\tau}}/V_A(\phi(\tCD_B)_{p_{\sigma\tau}}).$$
We have the following equalities:
\begin{enumerate}
\item $\dim_k (\tCD_A)_{p_{\tau}}/\phi(\tCD_B)_{p_{\tau}}$ is one if
$\tau$ is in $S$, and zero otherwise.
\item $\dim_k (\tCD_A)_{p_{\sigma\tau}}/\phi(\tCD_B)_{p_{\sigma\tau}}$
is one if $\sigma\tau$ is in $S$, and zero otherwise.
\item $\dim_k (\tCD_A)_{p_{\tau}}/V_A((\tCD_A)_{p_{\sigma\tau}}) = 1$.
\end{enumerate}
The desired formula for the rank of $\Lie(B/T)_{p_{\tau}}$ follows
immediately.

The Weil pairing induced by the polarization on $B$ identifies
$\Lie(B/T)^*_{q_{\tau}}$ with $(\Lie(B/T)^*)^{\perp}_{p_{\tau}},$
where $\perp$ denotes orthogonal complement with respect to the
Weil pairing
$$H^1_{\DR}(B/T)_{p_{\tau}} \times H^1_{\DR}(B/T)_{q_{\tau}} 
\rightarrow \OO_T.$$
Thus we conclude that $\Lie(B/T)_{q_{\tau}}$ has rank equal to
$2 - \rank \Lie(B/T)_{p_{\tau}}$, for all $\tau$.
\end{proof}

We are now in a position to show:
\begin{proposition} 
The map $\pi: Y_S \rightarrow X_S$ is an isomorphism.
\end{proposition}
\begin{proof}
This map is clearly proper, and $X_S$ is smooth, so it suffices to show that 
$\pi$ is bijective on points and injective on tangent spaces.

Let $x$ be a closed point of $X_S$, with residue field $k_x$.  Then
$k_x$ is perfect, and $x$ corresponds to a tuple $(A,\lambda,\rho)$
in $X_S(k_x)$.  We first construct a point of $Y_S$ lying over $x$.

Define a submodule $M$ of $H^1_{\DR}(A/k_x)$ by setting
$M_{p_{\tau}} = H^1_{\DR}(A/k_x)_{p_{\tau}}$ for $\tau$ not in $S$,
and $M_{p_{\tau}} = \Lie(A/k_x)_{p_{\tau}}^*$ for $\tau$ in $S$.  Set
$M_{q_{\tau}} = M_{p_{\tau}}^{\perp}$ for all $\tau$, where as usual
$\perp$ denotes orthogonal complement with respect to the Weil pairing:
$$H^1_{\DR}(A/k_x)_{p_{\tau}} \times H^1_{\DR}(A/k_x)_{q_{\tau}} \rightarrow 
k_x.$$

Let $\CD_A$ denote the Dieudonn{\'e} module of $A[p]$.  Then
$\CD_A$ is naturally isomorphic to $H^1_{\DR}(A/k_x)$.
We consider $M$ as a subspace of $\CD_A$ via this isomorphism.  It is
easy to verify that $M$ is stable under $F$ and $V$.

We have an exact sequence:
$$0 \rightarrow M \rightarrow \CD_A \rightarrow \CD_K \rightarrow 0,$$
where $\CD_K$ is the Dieudonn{\'e} module of some maximal isotropic
subgroup scheme $K$ of $A[p]$.  Let $B$ be the quotient $A/K$, and let $\phi$ be the natural
map of $A$ onto $B$.  Note that the submodule $\phi(H^1_{\DR}(B/k_x))$
of $H^1_{\DR}(A/k_x)$ is simply $M$.

Since $K$ is a maximal isotropic subgroup of $A[p]$,
it follows that there is a prime-to-$p$ polarization $\lambda^{\prime}$
on $B$ such that $p\lambda = \phi^{\vee}\lambda^{\prime}\phi$.
Thus $(A,\lambda,\rho,B,\lambda^{\prime},\phi)$ is a point of $Y_S$
mapping to $x$.

On the other hand, if $(A,\lambda,\rho,B, \lambda_{\prime}, \phi)$
is {\em any} point of $Y_S$ mapping to $x$, then by definition
$\phi(H^1_{\DR}(B/k_x))$ is equal to $M$.  Thus the kernel of
$\phi$ is equal to $K$, and hence this tuple is isomorphic
to the tuple constructed above.  In particular, $\pi$ is a bijection
on closed points.

It now suffices to show that the map $\pi$ is an injection on tangent 
spaces at every closed point.

Let $x$ be a closed point of $Y_S$, with (necessarily perfect)
residue field $k_x$.  The point $x$ corresponds to a tuple
$(A,\lambda,\rho,B,\lambda^{\prime},\phi)$. Let $\II$
denote the scheme $\spec k_x[\epsilon]/\epsilon^2$.  As 
in the proof of Proposition~\ref{prop:smooth}, specifying
a tangent vector to $X$ at $\pi(x)$ is equivalent to specifying, for each
$\tau$, a lift of $\Lie(A/k_x)^*_{p_{\tau}}$ to a subbundle 
$\omega_{A,\tau}$ of $(H^1_{\cris}(A/k_x)_{\II})_{p_{\tau}}.$

Fix choices of $\omega_{A,\tau}$ corresponding to a fixed tangent vector $v$
to $X_S$ at $x$, and let $\tA$ be the abelian scheme over
$\II$ corresponding to this tangent
vector. 
To lift this tangent vector to a tangent vector to $Y_S$ at
$x$, we must specify a lift $\tB$ of $B$, a lift of the polarization
$\lambda^{\prime}$ on $B$, and a lift of the isogeny $A \rightarrow B$. 

By Theorem~\ref{thm:grot}, to specify $\tB$, we must
specify for each $\tau$ lifts $\omega_{B,\tau}$
of $\Lie(B/k_x)^*_{p_{\tau}}$ and $\omega_{B,\tau}^{\prime}$
of $\Lie(B/k_x)^*_{q_{\tau}}$ to
$(H^1_{\cris}(B/k_x)_{\II})_{p_{\tau}}$ and
$(H^1_{\cris}(B/k_x)_{\II})_{q_{\tau}}$, respectively.  
The polarization $\lambda^{\prime}$ will lift to $\tB$
if and only if $\omega_{B,\tau}^{\prime} = \omega_{B,\tau}^{\perp}$ for
all $\tau$, where $\perp$ denotes orthogonal complement under
the pairing
$$
(H^1_{\cris}(B/k_x)_{\II})_{p_{\tau}} \times
(H^1_{\cris}(B/k_x)_{\II})_{q_{\tau}}
\rightarrow k_x[\epsilon]/\epsilon^2.
$$
Moreover, if a lift exists at all, it will be unique.

The map $\phi: A \rightarrow B$ will lift to a map $\tA \rightarrow \tB$
if and only if $\phi(\omega_{B,\tau}) \subset \omega_{A,\tau}$ 
and $\phi(\omega_{B,\tau}^{\prime}) \subset \omega_{A,\tau}^{\perp}$
for all $\tau$.  Again, if such a lift exists it will be unique.

Thus specifying a lift of $v$ to a tangent vector to 
$Y_S$ at $x$ is equivalent to specifying, for each $\tau$,
an $\omega_{B,\tau}$ such that $\phi(\omega_{B,\tau}) \subset \omega_{A,\tau}$.

Suppose first that $\tau$ is in $S$.  Then Lemma~\ref{lemma:lierank}
implies that the rank of $\Lie(B/k_x)^*_{p_{\tau}}$ is two.  In
particular, it is all of $H^1_{\DR}(B/k_x)_{p_{\tau}}$.
Thus the only possible choice for $\omega_{B,\tau}$ is all of
$(H^1_{\cris}(B/k_x)_{\II})_{p_{\tau}}$.

Suppose that $\sigma\tau$ is in $S$.  Then $\tau$ is not in $S$, as $S$
is sparse, and Lemma~\ref{lemma:lierank} implies that 
$\Lie(B/k_x)^*_{p_{\tau}} = 0$.  Thus the only possible choice
for $\omega_{B,\tau}$ is the zero module.

Finally, suppose that neither $\tau$ nor $\sigma\tau$ is in $S$.  Then
$\Lie(B/k_x)^*_{p_{\tau}}$ has rank one by Lemma~\ref{lemma:lierank}, so 
$\omega_{B,\tau}$ must have rank one, and map to
$\omega_{A,\tau}$ under $f$.  But as $\tau$ is not in $S$, $\phi$ induces
an isomorphism
$$
(H^1_{\cris}(B/k_x)_{\II})_{p_{\tau}} \rightarrow
(H^1_{\cris}(A/k_x)_{\II})_{p_{\tau}}.
$$
Moreover, $\omega_{A,\tau}$ has rank one.  Thus $\omega_{A,\tau}$ determines
$\omega_{B,\tau}$ uniquely.

It follows that the tangent space to $Y_S$ at $x$ injects
into the tangent space to $X_S$ at $\pi(x)$.  Hence $\pi$ is proper,
bijective on points, and injective on tangent spaces, and is 
therefore an isomorphism. 
\end{proof}

The isomorphism $Y_S \rightarrow X_S$ allows us to associate
to every triple $(A, \lambda, \rho) \in X_S$ 
a canonical isogeny $f: A \rightarrow B$, where $B$ is an
abelian scheme over $T$ with $\OO_F$ action and polarization 
$\lambda^{\prime}$.  The pair $(B,\lambda^{\prime})$ is ``almost'' a 
point of a $U(2)$ Shimura variety; the only thing missing is a
level structure $\rho^{\prime}$.  We cannot define such a level structure
in terms of $\CV$, because $r_{\tau}(\CV) = 1$ for all $\tau$, but, for
instance, $\Lie(B/T)_{p_{\tau}}$ has dimension $2$ if $\tau$ is in $S$.

Define $r_{\tau}^{\prime}$ and $s_{\tau}^{\prime}$ via:
\begin{enumerate}
\item $r_{\tau}^{\prime} = 2$ if $\tau \in S$.
\item $r_{\tau}^{\prime} = 0$ if $\sigma\tau \in S$.
\item $r_{\tau}^{\prime} = 1$ if neither $\tau$ nor $\sigma\tau$ lies in $S$.
\item $s_{\tau}^{\prime} = 2 - r_{\tau}^{\prime}$ for all $\tau$.
\end{enumerate}
Then, by Lemma~\ref{lemma:lierank}, $\Lie(B/T)_{p_{\tau}}$ has rank
$r_{\tau}^{\prime}$ for all $\tau$, and $\Lie(B/T)_{q_{\tau}}$ has
rank $s_{\tau}^{\prime}$ for all $\tau$.  Thus if we wish to make $B$ into
a point on some $U(2)$ Shimura variety $X_U(\CV^{\prime})$, we must have
$r_{\tau}(\CV^{\prime}) = r_{\tau}^{\prime}$ for all $\tau$. 

Fix an $F$-vector space $\CV^{\prime}$ with a pairing 
$\langle,\rangle^{\prime}$
satisfying the conditions of section~\ref{sec:basic} such that 
$r_{\tau}(\CV^{\prime}) = r_{\tau}^{\prime}$ for all $\tau$, and such
that there is an isomorphism
$\CV^{\prime} \otimes \AA^f_{\QQ} \cong \CV \otimes \AA^f_{\QQ}$
taking $\langle,\rangle^{\prime}$ to $\langle,\rangle$.
(The existence of such a vector space is a consequence of 
Corollary~\ref{cor:pairing}, which is proven in the appendix,
or of Proposition~\ref{prop:pairing} if we are working more generally
with $X_U^D(\CV)$ rather than $X_U(\CV)$.)
Fix in addition a particular isomorphism $\psi$ of 
$\CV^{\prime} \otimes \AA^f_{\QQ}$ with $\CV \otimes \AA^f_{\QQ}$, and
let $\CT^{\prime} = \psi(\CT)$.  Let $G^{\prime}$ denote the unitary
group associated to $\CV^{\prime}$ and $\langle,\rangle^{\prime}$.  The
map $\psi$ then defines an isomorphism $G(\AA^f_{\QQ}) \rightarrow
G^{\prime}(\AA^f_{\QQ})$, and we denote by $U^{\prime}$ the subgroup
of $G^{\prime}(\AA^f_{\QQ})$ that is the image of $U$ under this 
isomorphism.
Then $\psi$, the isogeny $f: A \rightarrow B$, and 
the $U$-level structure $\rho$ on $A$ induce a unique 
$U^{\prime}$-level structure $\rho^{\prime}$ on $B$.  In particular,
the triple $(B,\lambda^{\prime},\rho^{\prime})$ is a point
of $X^{\prime} = X_{U^{\prime}}(\CV^{\prime})_{k_0}$.  We thus
obtain a map $X_S \rightarrow X^{\prime}$.

The scheme $X_S$ has dimension $d-\#S$, whereas  we have:
$$\dim X^{\prime} = \sum_{\tau} r_{\tau}^{\prime}s_{\tau}^{\prime} = d - 2\#S.$$
Thus this map has nontrivial fibers, of dimension at least equal to the 
order of $S$.

In order to recover $(A,\lambda,\rho)$ given 
$(B,\lambda^{\prime},\rho^{\prime})$, we need to remember some extra data.
Let $(A,\lambda,\rho,B,\lambda^{\prime},\phi)$ be the point of $Y_S$
associated 
to $(A,\lambda,\rho)$.  Then the kernel of $\phi$ is contained in $A[p]$, so
multiplication by $p$ factors through $\phi$.  We thus obtain an isogeny
$\phi^{\prime}: B \rightarrow A$, of degree $p^{2d}$.

We make the following observations, which are clear::
\begin{itemize}
\item The kernel of $\phi^{\prime}$ is contained in $B[p]$.
\item We have
$p\lambda^{\prime} = (\phi^{\prime})^{\vee}\lambda \phi^{\prime}.$
\item $\phi^{\prime}(H^1_{\DR}(A/T)_{p_{\tau}})$ has
rank one for $\tau$ in $S$, and is zero for $\tau$ outside $S$.
\end{itemize}

Let $Y_S^{\prime}$ denote the moduli space that sends a $k_0$-scheme $T$
to the set of isomorphism classes of tuples $(B,\lambda^{\prime},\rho^{\prime},
A,\lambda, \phi^{\prime})$, where
\begin{enumerate}
\item $(B,\lambda^{\prime},\rho^{\prime})$ is a point of 
$X^{\prime}(T)$,
\item $A$ is an abelian scheme over $T$ of dimension $2d$, with an action
of $\OO_F$,
\item $\lambda$ is a polarization of $A$, of degree prime to $p$, and
\item $\phi^{\prime}: B \rightarrow A$ is an $\OO_F$-isogeny,
such that:
\begin{itemize}
\item $p\lambda^{\prime} = (\phi^{\prime})^{\vee}\lambda\phi^{\prime}$, and
\item for each $\tau$, $\phi^{\prime}(H^1_{\DR}(A/T)_{p_{\tau}})$ is
locally free of rank one if $\tau$ is in $S$, and zero otherwise.
\end{itemize}
\end{enumerate}

The construction outlined above defines a map $Y_S \rightarrow
Y_S^{\prime}$. Let us now define an inverse map.
Suppose we have a point $(B,\lambda^{\prime},\rho^{\prime},A,\lambda,\phi^{\prime})$
of $Y_S^{\prime}$.  As before, the $U^{\prime}$-level structure on
$B$ induces a unique $U$-level structure $\rho$ on $A$, and hence
a unique point $(A,\lambda,\rho)$ of $X$.  Factoring
multiplication by $p$ on $B$ through $\phi^{\prime}$ gives us an isogeny
$\phi: A \rightarrow B$, of degree $p^{2d}$. 

\begin{lemma}
The tuple $(A,\lambda,\rho,B,\lambda^{\prime},\phi)$ defines a point
of $Y_S$ that maps to $(B,\lambda^{\prime},p\rho^{\prime},A,\lambda,\phi^{\prime})$
under the map $Y_S \rightarrow Y_S^{\prime}$.  Thus $Y_S$ is isomorphic to $Y_S^{\prime}$. 
\end{lemma}
\begin{proof}
Exactly as above, we see that the kernel of $\phi$ is contained in $A[p]$,
that $p\lambda = \phi^{\vee}\lambda^{\prime} \phi$, and that
$\phi(H^1_{\DR}(B/T)_{p_{\tau}})$ has rank one if $\tau$ is in $S$ and
rank two otherwise.  It remains to show that for $\tau \in S$, we have
$$\phi(H^1_{\DR}(B/T)_{p_{\tau}}) = \Lie(A/T)^*_{p_{\tau}},$$
and also 
$$\phi(H^1_{\DR}(B/T)_{p_{\tau}}) = \Fr(H^1_{\DR}(A^{(p)}/T)_{p_{\tau}}).$$

Note that for $\tau$ in $S$,
$\Lie(A/T)^*_{p_{\tau}}$ has rank one, $\Lie(B/T)^*_{p_{\tau}}$ has
rank two, and the kernel of the map 
$$H^1_{\DR}(B/T)_{p_{\tau}} \rightarrow H^1_{\DR}(A/T)_{p_{\tau}}$$ 
has rank one.  Thus $\phi$ induces a surjection of 
$\Lie(B/T)^*_{p_{\tau}}$ (equivalently, of $H^1_{\DR}(B/T)_{p_{\tau}}$) 
onto $\Lie(A/T)^*_{p_{\tau}}$.

On the other hand, as $\phi$ commutes with relative Frobenius, we have
$$\phi(\Fr(H^1_{\DR}(B^{(p)}/T)_{p_{\tau}})) \subset
\Fr(H^1_{\DR}(A^{(p)}/T)_{p_{\tau}}).$$
For $\tau \in S$, $\Fr(H^1_{\DR}(B^{(p)}/T)_{p_{\tau}})$ has
rank two, $\Fr(H^1_{\DR}(A^{(p)}/T)_{p_{\tau}})$ has rank one,
and the kernel of $\phi$ has rank one, so again $\phi$ induces a surjection
of $H^1_{\DR}(B/T)_{p_{\tau}}$ onto $\Fr(H^1_{\DR}(A^{(p)}/T)_{p_{\tau}}).$

This construction defines a map $Y_S^{\prime} \rightarrow Y_S$
that is nearly an inverse to the map $Y_S \rightarrow Y_S^{\prime}$
defined earlier.  (The composition of the two maps is the automorphism of
$Y_S$ which multiplies the level structure $\rho$ by $p$.)
We thus obtain an isomorphism $Y_S \rightarrow Y_S^{\prime}.$
\end{proof}

Let $Z_S$ be the moduli space of tuples
$(B,\lambda^{\prime},\rho^{\prime},M)$, where 
$(B,\lambda^{\prime},\rho^{\prime})$ is a point of $X^{\prime}(T)$,
and $M$ is an $\OO_F$-stable subbundle of $H^1_{\DR}(B/T)$ such that:
\begin{enumerate}
\item $M_{p_{\tau}}$ is locally free of rank one for $\tau$ in $S$, and zero
for $\tau$ not in $S$, and
\item $M_{q_{\tau}} = M_{p_{\tau}}^{\perp}$, where $\perp$ is with respect to
the Weil pairing
$$H^1_{\DR}(B/T)_{p_{\tau}} \times H^1_{\DR}(B/T)_{q_{\tau}} \rightarrow \OO_T.$$
\end{enumerate}

Let $\CB$ denote the universal abelian scheme on $X^{\prime}$.
Then $Z_S$ is simply the fiber product, over $X^{\prime}$,
of the bundles $P(H^1_{\DR}(\CB/X^{\prime})_{p_{\tau}})$, for $\tau \in S$.

We have a natural map $Y_S^{\prime} \rightarrow Z_S$ that
takes the point
$(B,\lambda^{\prime},\rho^{\prime},A,\lambda,\phi^{\prime})$ of
$Y_S^{\prime}$ to the point
$(B,\lambda^{\prime},\rho^{\prime},\phi^{\prime}(H^1_{\DR}(A/T)))$ of
$Z_S$.

\begin{proposition} \label{prop:pointstratum}
The natural map $Y_S^{\prime} \rightarrow Z_S$ is
a bijection on points.
\end{proposition}
\begin{proof} We construct an inverse map on the level of points.  Let $k$
be a perfect field containing $k_0$, and let 
$(B,\lambda^{\prime},\rho^{\prime},M)$ be a $k$-valued point of $Z_S$.

Let $\CD_B$ denote the Dieudonn{\'e} module of $B[p]$, and consider
$M$ as a subspace of $\CD_B$.  Then $M$ is stable under $F$ and $V$,
as for $\tau$ in $S$, $F((\CD_B)_{p_{\tau}})$ and $V((\CD_B)_{p_{\tau}})$
both vanish.  Hence we have an exact sequence of Dieudonn{\'e} modules:
$$0 \rightarrow M \rightarrow \CD_B \rightarrow \CD_K \rightarrow 0,$$
where $\CD_K$ is the Dieudonn{\'e} module of some finite flat group scheme
$K$.  We thus obtain an injection of $K$ into $B[p]$.

Let $A = B/K$, and let $\phi^{\prime}: B \rightarrow A$ be the quotient
map.  Since $K$ is a maximal isotropic subgroup of $B[p]$, there
is a prime-to-$p$ polarization $\lambda$ on $A$ such that
$p\lambda^{\prime} = (\phi^{\prime})^{\vee}\lambda\phi^{\prime}.$
Then $(B,\lambda^{\prime},\rho^{\prime},A,\lambda,\phi^{\prime})$ is
a point of $Y_S^{\prime}$, and it is easy to see that this gives
an inverse to the map $Y_S^{\prime} \rightarrow Z_S$ defined above.
\end{proof}

Composing the map $Y_S^{\prime} \rightarrow Z_S$ with the
sequence of isomorphisms $X_S \cong Y_S \cong Y_S^{\prime}$ yields
a map $X_S \rightarrow Z_S$ that is a bijection on points.  (In
general this map fails to be a geometric isomorphism.)

\begin{definition} A map $Y \rightarrow Z$ of $k$-schemes is a {\em Frobenius factor}
if there exists a map $f^{\prime}: Z_{p^r} \rightarrow Y$ such
that $ff^{\prime}: Z_{p^r} \rightarrow Z$ is the $r$th power of
the Frobenius morphism.  (Here $Z_{p^r}$ is the $k$-scheme obtained
from $Z$ by composing the structure map $Z \rightarrow \spec k$
with the $p^r$-th power Frobenius automorphism of $\spec k$.) 
\end{definition}

Note that Frobenius factors are isomorphisms on {\'e}tale cohomology.

\begin{proposition} \label{prop:bijection}
Let $Y$ and $Z$ be schemes of finite type over a perfect
field $k$ of characteristic $p$, such that $Z$ is normal and $Y$ is reduced.  
Let 
$f: Y \rightarrow Z$ be a proper morphism that is a bijection on points.
Then $f$ is a Frobenius factor.
\end{proposition}
\begin{proof}
As $f$ is a bijection on points, $K(Y)$ is a purely inseparable extension
of $K(Z)$.  Thus
for $r$ sufficiently large, every $p^r$-th power of an element of
$K(Y)$ lies in $K(Z)$.  In particular, we have $K(Y) \subset K(Z_{p^r}),$
For each $U \subset Z$ in an affine cover,
$\OO_{Z_{p^r}}(U)$ is the integral closure of $\OO_Z(U)$ in $K(Z_{p^r})$.
Since $f: Y \rightarrow Z$ is finite, $\OO_Y(U)$ is integral over $\OO_Z(U)$,
and hence contained in $\OO_{Z_{p^r}}(U)$.  These inclusions give the
map $f^{\prime}: Z_{p^r} \rightarrow Y$, as required.
\end{proof}

Combining this with the preceding result, we obtain:
\begin{theorem} \label{thm:stratum}
The natural map $X_S \rightarrow Z_S$ is a Frobenius factor.
\end{theorem}

We will have need of the following lemma in the next section:

\begin{lemma} \label{lemma:transfer}
Let $\CA_S$ be the universal abelian scheme on $X_S$.  Then for all
$\tau \notin S$, the bundle $H^1_{\DR}(\CA_S/X_S)_{p_{\tau}}$
is the pullback of $H^1_{\DR}(\CB/X^{\prime})_{p_{\tau}}$ under the
natural map $X_S \rightarrow X^{\prime}$ defined above.
\end{lemma}
\begin{proof}
Let $(\CA_S,\lambda,\rho, \CB_S,\lambda^{\prime},\phi_S)$ be the universal
object of $Y_S$. 
Then $\phi_S$ induces an isomorphism:
$$\phi_S: H^1_{\DR}(\CA_S/X_S)_{p_{\tau}} \cong 
H^1_{\DR}(\CB_S/X_S)_{p_{\tau}}$$
for all $\tau$ outside $S$.

On the other hand, by the definition of the natural map 
$X_S \rightarrow X^{\prime}$, we have 
$\CB_S \cong \CB \times_{X^{\prime}} X_S$.  Thus 
$H^1_{\DR}(\CB_S/X_S)_{p_{\tau}}$ is naturally isomorphic to the pullback
of $H^1_{\DR}(\CB/X^{\prime})_{p_{\tau}}$ for all $\tau$.
\end{proof}

\section{$\Gamma_0(p)$ level structures} \label{sec:DR}

We now turn to the reduction of $U(2)$ Shimura varieties with level
structure at $p$.  In particular fix a $\CV$ such that $r_{\tau}(\CV) =
s_{\tau(\CV)} = 1$ for all $\tau$, and consider the Shimura variety
$X_{U,p}(\CV)$ defined in section~\ref{sec:basic}.  Let
$\CA$ and $\CB$ be the universal abelian schemes on $X_{U,p}(\CV)$,
and let $\phi_{\univ}: \CA \rightarrow \CB$ denote the universal isogeny
between them.  Multiplication by $p$ on $\CA$ factors through $\phi_{\univ}$, 
and we denote by $\phi^{\prime}_{\univ}$ the isogeny from $\CB$ to $\CA$ 
obtained in this manner.

Let $\CE_{p_{\tau}}$ and $\CE^{\prime}_{p_{\tau}}$ be the bundles
$\Lie(\CA/X_{U,p}(\CV))^*_{p_{\tau}}$ and 
$\Lie(\CB/X_{U,p}(\CV))^*_{p_{\tau}}$, respectively. 
For all $\tau$, $\phi_{\univ}$ and $\phi^{\prime}_{\univ}$
can be thought of as sections of
$\Hom(\CE^{\prime}_{p_{\tau}}, \CE_{p_{\tau}})$
and $\Hom(\CE_{p_{\tau}}, \CE^{\prime}_{p_{\tau}})$,
which we denote by $h_{\tau}$ and $h_{\tau}^{\prime}$, respectively.

Note that the composition of $h_{\tau}$ and $h_{\tau}^{\prime}$ in either
direction is multiplication by $p$, because the same holds for $\phi_{\univ}$ and 
$\phi_{\univ}^{\prime}$.  Thus in particular on any point $x$ of the fiber of 
$X_{U,p}(\CV)$ over $p$, either $h_{\tau}$ or $h_{\tau}^{\prime}$
vanishes.

For subsets $S$ and $S^{\prime}$ of $\CS$, define the
subscheme $X_{U,p}(\CV)_{S,S^{\prime}}$ of $X_{U,p}(\CV)_{k_0}$  
to be the common zeros of the $h_{\tau}$ for $\tau \in S$ and
the $h_{\tau}^{\prime}$ for $\tau \in S^{\prime}$.
Then $X_{U,p}(\CV)_{k_0}$ is the union of the 
$X_{U,p}(\CV)_{S,S^{\prime}}$ for pairs of sets $S,S^{\prime}$
with $S \cup S^{\prime} = \CS$.

We study these subschemes by constructing a local
model for $X_{U,p}(\CV)$ in which the above subschemes have natural
analogues.  Constructing such models is standard; the technique is due originally
to de Jong~\cite{sgee}.  The result we need is a special case of results of 
Rapoport-Zink~\cite{RZ1}, Appendix.

For each $\tau$, let $M_{\tau}$ and $M_{\tau}^{\prime}$ be free $W(\FF_p)$ 
modules of rank two, with bases 
$\{e_{1,\tau}$, $e_{2,\tau}\}$ and $\{e_{1,\tau}^{\prime}$, $e_{2,\tau}^{\prime}\}$, 
respectively.  Let
$\alpha_{\tau}: M_{\tau} \rightarrow M_{\tau}^{\prime}$ be defined by 
$\alpha_{\tau}(e_{1,\tau}) = e_{1,\tau}^{\prime}$
and $\alpha_{\tau}(e_{2,\tau}) = pe_{2,\tau}^{\prime}$.  Similarly, define 
$\alpha^{\prime}_{\tau}$
by $\alpha^{\prime}_{\tau}(e_{1,\tau}^{\prime}) = pe_{1,\tau}$ and 
$\alpha^{\prime}_{\tau}(e_{2,\tau}^{\prime}) = e_{2,\tau}$.

Define $\tX_{U,p}(\CV)$ to be the scheme representing the functor
parameterizing tuples 
$(A,\lambda,\rho,B,\lambda^{\prime},\phi,\{e_\tau\},
\{e_{\tau}^{\prime}\}),$ where
\begin{enumerate}
\item $(A,\lambda,\rho,B,\lambda^{\prime},\phi)$ is a point of 
$X_{U,p}(\CV)(T)$,
\item for each $\tau$, $e_{\tau}$ and $e_{\tau}^{\prime}$ are isomorphisms of
$H^1_{\DR}(B/T)_{p_{\tau}}$ with $M_{\tau} \otimes_{W(k_0)} T$ and $H^1_{\DR}(A/T)_{p_{\tau}}$ 
with $M_{\tau}^{\prime} \otimes_{W(k_0)} T$, respectively.
\item We require that under the isomorphism $e_{\tau}$ and $e_{\tau}^{\prime}$,
the map 
$$\phi: H^1_{\DR}(B/T)_{p_{\tau}} \rightarrow H^1_{\DR}(A/T)_{p_{\tau}}$$
corresponds to $\alpha$, and the map 
$$\phi^{\prime}: H^1_{\DR}(A/T)_{p_{\tau}} \rightarrow H^1_{\DR}(B/T)_{p_{\tau}}$$
corresponds to $\alpha^{\prime}$.  (Here $\phi^{\prime}: B \rightarrow A$ is
the map obtained by factoring multiplication by $p$ through $\phi$.)
\end{enumerate}

It is then standard that 
the natural map to $\tX_{U,p}(\CV) \rightarrow X_{U,p}(\CV)$ is smooth and surjective.

We now describe the local model for $X_{U,p}$.  Let $\CM$
denote the $W(k_0)$-scheme parameterizing tuples
$(\{V_{\tau},V^{\prime}_{\tau}\})$, where for each $\tau$,
$V_{\tau}$ is a rank one subbundle of $(M_{\tau})_T$, $V^{\prime}_{\tau}$ is a rank one 
subbundle of $(M_{\tau}^{\prime})_T$, $\alpha_{\tau}$ maps $V$ to $V^{\prime}$, and 
$\alpha^{\prime}_{\tau}$ maps $V^{\prime}$ to $V$.  The scheme $\CM$ is the product of
schemes $\CM_{\tau}$ paramaterizing $(V_{\tau},V_{\tau}^{\prime})$ as above for a single
$\tau$.

The general fiber of $\CM_{\tau}$ is isomorphic to $\PP^1$, while the fiber of
$\CM_{\tau}$ over $\FF_p$ is the union of two components.  More precisely,
let $\beta_{\tau}$, $\beta^{\prime}_{\tau}$ be the maps 
$V_{\tau} \rightarrow V_{\tau}^{\prime}$ and $V_{\tau}^{\prime} \rightarrow V_{\tau}$
induced by $\alpha_{\tau}$ and $\alpha_{\tau}^{\prime}$, respectively.
These are sections of tautological bundles on $\CM_{\tau}$.
Their composition is multiplication by $p$, and the special fiber of
$\CM_{\tau}$ is thus the union of $Z(\alpha_{\tau})$ and $Z(\beta_{\tau})$.
Each of these components is isomorphic to $\PP^1_{\FF_p}$.

For $S$ and $S^{\prime}$ such that $S \cup S^{\prime} = \CS$,
define $\CM_{S,S^{\prime}}$ to be the zero locus of $\beta_{\tau}$
for $\tau \in S$ and $\beta_{\tau}^{\prime}$ for $\tau \in S^{\prime}$.

We have a map $g: \tX_{U,p}(\CV) \rightarrow \CM$ 
that takes $(A,\lambda,\rho,B,\lambda^{\prime},\phi,\{e_{\tau}\},\{e_{\tau}^{\prime}\})$
to $(\{e_{\tau}(\Lie(B/T)^*_{p_{\tau}}),e_{\tau}^{\prime}(\Lie(A/T)^*_{p_{\tau}})\})$.
The arguments of~\cite{RZ1} show that $g$ is smooth, 
and hence (via~\cite{sgee}, corollary 4.6) that
every point $x$ of $X_{U,p}(\CV)$ has an {\'e}tale neighborhood that 
is isomorphic to an {\'e}tale neighborhood of some point $w$ of $\CM$.  
Moreover, the pullbacks to $\tX_{U,p}(\CV)$ of $\beta_{\tau}$ and $\beta^{\prime}_{\tau}$
coincide with the pullbacks of $h_{\tau}$ and $h_{\tau}^{\prime}$.  Thus
if $x$ lies in $X_{U,p}(\CV)_{S,S^{\prime}}$, then $w$ can be
taken to be inside $\CM_{S,S^{\prime}}$.  The following corollaries are
immediate.

\begin{corollary}
$X_{U,p}(\CV)_{k_0}$ is reduced.
\end{corollary}

\begin{corollary} 
Let $S,S^{\prime} \subset \CS$ satisfy 
$S \cup S^{\prime} = \CS$.  Then
$X_{U,p}(\CV)_{S,S^{\prime}}$ is a smooth subscheme of
$X_{U,p}(\CV)_{k_0}$, of codimension equal
to $\# (S \cap S^{\prime})$.  In particular, if $S$ and $S^{\prime}$
are disjoint then $X_{U,p}(\CV)_{S,S^{\prime}}$ is a disjoint
union of irreducible components of $X_{U,p}(\CV)$.  Moreover,
every irreducible component of $(X_{U,p})_{k_0}$ is smooth.
\end{corollary}

\begin{corollary} 
Let $x$ be a closed point of the fiber of 
$X_{U,p}(\CV)$.  Let $S$ be the set of all $\tau$ such that
$h_{\tau}(x) = 0$, and let $S^{\prime}$ be the set of all $\tau$
such that $h_{\tau}^{\prime}(x) = 0$.  Then the completion of the strict
henselization of the local ring of $X_{U,p}(\CV)$ at $x$ is
isomorphic to:
$$W(\overline{\FF}_p)[[X_{\tau}: \tau \in S, Y_{\tau}: \tau \in S^{\prime}]]/
(X_{\tau}Y_{\tau} - p: \tau \in S \cap S^{\prime}).$$
\end{corollary}

\begin{remark} \rm
In visualizing the way in which the subvarieties
$X_{U,p}(\CV)_{S,S^{\prime}}$ fit together, it is helpful to
consider them as corresponding to various faces of the cube
$[0,1]^{\CS}$.  Let 
$X_{U,p}(\CV)_{S,S^{\prime}}$ correspond to the face consisting
of points whose $\tau$ coordinate is zero for all $\tau$ in $S$ but not in
$S^{\prime}$, and one for all $\tau$ in $S^{\prime}$ but not in $S$.
Then the dimension of any face of this cube is equal to the {\em co}dimension
of the corresponding subvariety.  Moreover, if $X_1, \dots, X_r$ is a 
collection of subschemes of the above form, corresponding to faces
$F_1, \dots, F_r$ of the cube, then
the intersection of the $X_i$ is also of the above form, and corresponds to
the smallest face of $[0,1]^{\CS}$ containing all of the 
$F_i$. 
\end{remark}

\begin{remark} \rm
Strictly speaking, we have only proven that 
$X_{U,p}(\CV)_{S,S^{\prime}}$ has the above properties {\em if
it is nonempty.}  This nonemptiness is a
consequence of our characterization of the closed points on
$X_{U,p}(\CV)_{S,S^{\prime}}$, below.
\end{remark}

Now that we have an understanding of the local structure of the special
fiber, we turn to the global structure of the subvarieties
$X_{U,p}(\CV)_{S,S^{\prime}}$.  The key is the following
observation:

\begin{proposition}
Let $k$ be a perfect field of characteristic $p$, and $x$ a $k$-valued
point of $X_{U,p}(\CV)$, corresponding to an isogeny
$\phi_x: A/k \rightarrow B/k$.  Then:
\begin{enumerate}
\item $h_{\tau}^{\prime}(x) = 0$ if and only if 
$\phi_x(H^1_{\DR}(B/k))_{p_{\tau}} = \Lie(A/k)^*_{p_{\tau}}$, and
\item $h_{\tau}(x) = 0$ if and only if $\phi_x(H^1_{\DR}(B/k))_{p_{\sigma\tau}}
= \Fr(H^1_{\DR}(A^{(p)}/k))_{p_{\sigma\tau}}$, where $\Fr$ is the relative
Frobenius $A \rightarrow A^{(p)}$. 
\end{enumerate}
\end{proposition}
\begin{proof}
Let $\phi_x^{\prime}$ be the map $B/k \rightarrow A/k$ obtained by factoring
multiplication by $p$ on $A$ through $\phi_x$.
The condition that $h_{\tau}^{\prime}(x) = 0$ is equivalent to requiring
that the map 
$\phi_x^{\prime}: \Lie(A/k)^*_{p_{\tau}} \rightarrow \Lie(B/k)^*_{p_{\tau}}$
be zero.  This in turn is equivalent to requiring that the kernel of
$\phi_x^{\prime}$ on $H^1_{\DR}(A/k)_{p_{\tau}}$ contain 
$\Lie(A/k)^*_{p_{\tau}}$;
as both have dimension one they must be equal.  On the other hand, the
kernel of $\phi_x^{\prime}$ on $H^1_{\DR}(A/k)_{p_{\tau}}$ is equal to
$\phi_x(H^1_{\DR}(B/k))_{p_{\tau}}$, so $h_{\tau}^{\prime}(x) = 0$ 
if and only if
$\phi_x(H^1_{\DR}(B/k))_{p_{\tau}}$ is equal to $\Lie(A/k)^*_{p_{\tau}}$ as 
claimed.  

For the second claim, we use Dieudonn{\'e} theory.  Let $\tCD_A$ and
$\tCD_B$ denote the contravariant Dieudonn{\'e} modules of $A[p^{\infty}]$
and $B[p^{\infty}]$, respectively.  As usual, there is a natural isomorphism
$\tCD_A/p\tCD_A \cong H^1_{\DR}(A/k)$, and a similar isomorphism for $B$. 

The map $\phi_x$ induces an injection of $\tCD_B$ into $\tCD_A$, whose image
contains $p\tCD_A$; under the above
isomorphisms the map $\phi_x$ induces on de Rham cohomology is simply the 
mod $p$ reduction of this injection.  

Now by the proof of the first part of the proposition (and the symmetry
of $\phi_x$ and $\phi_x^{\prime})$, $h_{\tau}(x)$ is zero if, and only if,
$\phi_x^{\prime}(H^1_{\DR}(A/k)_{p_{\tau}})$ is equal to 
$\Lie(B/k)_{p_{\tau}}^*$.
In the language of Dieudonn{\'e} theory, this occurs if and only if 
$\phi_x^{\prime}(\tCD_A)_{p_{\tau}}$ is equal to
$V_B((\tCD_B)_{p_{\sigma\tau}})$.  Since $\phi_x$ is injective on Dieudonn{\'e}
modules, we can check this
equality after composition with $\phi_x$; i.e. $h_{\tau}(x)$ is zero if and only
if $p(\tCD_A)_{p_{\tau}}$ equals $\phi_x(V_B((\tCD_B)_{p_{\sigma\tau}}))$
(or, equivalently, if and only if $p(\tCD_A)_{p_{\tau}}$ equals
$V_A(\phi_x(\tCD_B)_{p_{\sigma\tau}})$). 

Finally, as $V_A$ is injective and $V_AF_A = p$, we have that
$h_{\tau}(x)$ is zero if and only if $\phi_x(\tCD_B)_{p_{\sigma\tau}}$
is equal to $F_A((\tCD_A)_{p_{\tau}})$.  In terms of de Rham cohomology,
this occurs if and only if $\phi_x(H^1_{\DR}(B/k)_{p_{\sigma\tau}})$ is
equal to $\Fr(H^1_{\DR}(A^{(p)}/k)_{p_{\sigma\tau}}).$
\end{proof}

In light of this proposition, we define $Y_{S,S^{\prime}}$ to be the
subscheme of $X_{U,p}(\CV)$ whose $T$-points are tuples
$(A,\lambda,\rho,B,\lambda^{\prime},\phi)$ that satisfy:
\begin{enumerate}
\item $\phi(H^1_{\DR}(B/T)_{p_{\tau}}) = \Lie(A/T)^*_{p_{\tau}}$ if
$\tau$ is in $S^{\prime}$, and
\item $\phi(H^1_{\DR}(B/T)_{p_{\tau}}) = \Fr(H^1_{\DR}(A^{(p)}/T)_{p_{\tau}})$
if $\sigma^{-1}\tau$ is in $S$.
\end{enumerate}

The proposition shows that the underlying point sets of 
$X_{U,p}(\CV)_{S,S^{\prime}}$ and $Y_{S,S^{\prime}}$
are the same.  Since $X_{U,p}(\CV)_{S,S^{\prime}}$ is
reduced, it follows that we have an isomorphism:
$$X_{U,p}(\CV)_{S,S^{\prime}} \cong  
Y_{S,S^{\prime}}^{\red}.$$

Now let $Z_{S,S^{\prime}}$ be the scheme representing the functor that
takes a $k_0$-scheme $T$ to the set of isomorphism classes of
tuples $(A,\lambda,\rho,M)$, where $(A,\lambda,\rho)$ is a point of
$X_U(\CV)$, and $M$ is an $\OO_F$-stable subbundle of $H^1_{\DR}(A/T)$ 
that satisfies:
\begin{enumerate}
\item $M_{p_{\tau}}$ is locally free of rank one for all $\tau$,
\item $M_{q_{\tau}} = M_{p_{\tau}}^{\perp}$ for all $\tau$,
\item $M_{p_{\tau}} = \Lie(A/T)_{p_{\tau}}^*$ if $\tau$ lies in $S^{\prime}$,
and
\item $M_{p_{\tau}} = \Fr(H^1_{\DR}(A^{(p)}/T)_{p_{\tau}})$ if 
$\sigma^{-1}\tau$ lies in $S$.
\end{enumerate}

The structure of $Z_{S,S^{\prime}}$ is straightforward to describe
explicitly.  It is a fiber product of $\PP^1$-bundles over a closed
stratum of $X_U(\CV)$.  Let $\Sigma_{S,S^{\prime}}$ be the set of all
$\tau$ with $\tau$ in $S^{\prime}$ and $\sigma^{-1}\tau$ in $S$.
If $(A,\lambda,\rho,M)$ is
a point of $Z_{S,S^{\prime}}$, then $(A,\lambda,\rho)$ lies in
the stratum $X_U(\CV)_{\Sigma_{S,S^{\prime}}}$.  

In fact, $Z_{S,S^{\prime}}$ is simply a product of $\PP^1$-bundles over
this stratum.
Let $\overline{\CA}$ be the universal abelian scheme
on $X_U(\CV)_{\Sigma_{S,S^{\prime}}}$, and for each $\tau$,
let $\overline{\CE}_{\tau}$ denote the rank two bundle 
$H^1_{\DR}(\overline{\CA}/
X_U(\CV)_{\Sigma_{S,S^{\prime}}})_{p_{\tau}}$.
Let $\Theta_{S,S^{\prime}}$ denote the set of all $\tau$ with $\tau$ {\em not}
in $S^{\prime}$ and $\sigma^{-1}\tau$ {\em not} in $S$.
Then $Z_{S,S^{\prime}}$ is simply the fiber product (over $
X_U(\CV)_{\Sigma_{S,S^{\prime}}}$) of the bundles
$\PP(\overline{\CE}_{\tau})$, for $\tau$ in $\Theta_{S,S^{\prime}}$.

There is an obvious map
$Y_{S,S^{\prime}} \rightarrow Z_{S,S^{\prime}}$ that takes
$(A,\lambda,\rho,B,\lambda^{\prime},\phi)$ to
$(A,\lambda,\rho,\phi(H^1_{\DR}(B/T)))$.

\begin{proposition} If $S \cup S^{\prime}$ is all of
$\CS$, then the map $Y_{S,S^{\prime}} \rightarrow
Z_{S,S^{\prime}}$ is a bijection on points.
\end{proposition}
\begin{proof}
We construct an inverse map.  Given a $k$-valued point
$(A,\lambda,\rho,M)$, of $Z_{S,S^{\prime}}$, with $k$ perfect, we
identify $M$ with a subspace of the Dieudonn{\'e} module $\CD_A$ of
$A[p]$.  Note that if $\tau$ is in $S^{\prime}$, then $M_{p_{\tau}}$
is equal to $V_A((\CD_A)_{p_{\sigma\tau}})$, and if $\sigma^{-1}\tau$ is
in $S$, then $M_{p_{\tau}}$ is equal to $F_A((\CD_A)_{p_{\sigma^{-1}\tau}})$.

For a given $\tau$, either $\sigma^{-1}\tau$ is in $S^{\prime}$ or
$\sigma^{-1}\tau$ is in $S$.  In the former case $M_{p_{\sigma^{-1}\tau}}$
is equal to $V_A((\CD_A)_{p_{\tau}})$, and hence clearly contains
$V_A(M_{p_{\tau}})$.  In the latter case 
$M_{p_{\tau}} = F_A((\CD_A)_{p_{\sigma^{-1}\tau}})$, and hence 
$V_A(M_{p_{\tau}}) = 0$.  Either way, we have $V_A(M_{p_{\tau}}) \subset
M_{p_{\sigma^{-1}\tau}}$.  It follows that $M$ is stable under $V_A$;
a similar argument shows it is stable under $F_A$.

Thus, exactly as in the proof of Proposition~\ref{prop:pointstratum},
we obtain a subgroup $K$ of $A$.  Let $B=A/K$, and let $\phi$ be the quotient
map $A \rightarrow B$.  We have a polarization $\lambda^{\prime}$
on $B = A/K$, such that $p\lambda = \phi^{\vee}\lambda^{\prime}\phi$,
and the point $(A,\lambda,\rho,B,\lambda^{\prime},\phi)$ is a point
of $Y_{S,S^{\prime}}$.  This map is easily checked to be an inverse map
(at the level of points) to the natural map $Y_{S,S^{\prime}} \rightarrow
Z_{S,S^{\prime}}$.
\end{proof}

\begin{theorem} \label{thm:DR1}
The natural morphism 
$X_{U,p}(\CV)_{S,S^{\prime}} \rightarrow Z_{S,S^{\prime}}$
is a Frobenius factor.
\end{theorem}
\begin{proof}
This is immediate from Proposition~\ref{prop:bijection} and the above
result.
\end{proof}

If $S$ and $S^{\prime}$ are disjoint, so that 
$X_{U,p}(\CV)_{S,S^{\prime}}$ has codimension zero in the special
fiber, then we can go further.

\begin{lemma}
If $S$ and $S^{\prime}$ are disjoint, then $\Sigma_{S,S^{\prime}}$ is
sparse.
\end{lemma}
\begin{proof}
If both $\tau$ and $\sigma\tau$ are in $\Sigma_{S,S^{\prime}}$, then
$\tau$ must be in both $S$ and $S^{\prime}$.
\end{proof}

In this case we can apply the results of the previous section.  Let
$\CV^{\prime}$ be a two-dimensional $F$-vector space with pairing such
that $\CV^{\prime}(\AA^f_{\QQ}) \cong \CV(\AA^f_{\QQ})$ but
for which
\begin{enumerate}
\item $r_{\tau}(\CV^{\prime}) = 2$ if $\tau \in \Sigma_{S,S^{\prime}}$,
\item $r_{\tau}(\CV^{\prime}) = 0$ if $\sigma\tau \in \Sigma_{S,S^{\prime}}$,
and
\item $r_{\tau}(\CV^{\prime}) = 1$ otherwise.
\end{enumerate}

Then, by Theorem~\ref{thm:stratum},  
$X_U(\CV)_{S_*(S,S^{\prime})}$ is (up to a Frobenius factor) a product
of $\PP^1$-bundles
over $X_U(\CV^{\prime})$.  Since $Z_{S,S^{\prime}}$ is (up to a Frobenius
factor) a product of $\PP^1$-bundles over $X_U(\CV)_{S_*(S,S^{\prime})}$,
we can express $Z_{S,S^{\prime}}$ as a product of $\PP^1$-bundles
over $X_U(\CV)^{\prime}$.  Lemma~\ref{lemma:transfer},
allows us to do this entirely in terms of natural bundles on $X_U(\CV^{\prime})$.  

More precisely, let $\overline{\CA}^{\prime}$ be the universal
abelian variety on $X_U(\CV^{\prime})$, and for each $\tau$ let
$\overline{\CE}^{\prime}_{\tau}$ be the vector bundle 
$H^1_{\DR}(\overline{\CA}^{\prime}/X_U(\CV^{\prime}))_{p_{\tau}}.$ 
Let $Z^{\prime}_{S,S^{\prime}}$ be the fiber product 
(over $X_U(\CV^{\prime})$) of the bundles $\PP(\CE^{\prime}_{\tau})$,
for $\tau$ in $\Sigma_{S,S^{\prime}} \cup \Theta_{S,S^{\prime}}$.
By Lemma~\ref{lemma:transfer} and Theorem~\ref{thm:stratum}
we obtain a map $Z_{S,S^{\prime}} \rightarrow Z^{\prime}_{S,S^{\prime}}$
that is a Frobenius factor.  Composing this with the map
$X_{U,p}(\CV)_{S,S^{\prime}} \rightarrow Z_{S,S^{\prime}}$, we obtain:

\begin{theorem} \label{thm:DR2}
The map $X_{U,p}(\CV)_{S,S^{\prime}}
\rightarrow Z^{\prime}_{S,S^{\prime}}$ described above is a Frobenius factor.
\end{theorem}

This gives a complete global description of the irreducible components
of $X_{U,p}(\CV)$ in terms of Shimura varieties for various inner forms of $G$.

\begin{example} \label{ex:quadratic} \rm
Suppose $F^+$ is real quadratic, and let $\tau_1$ and $\tau_2$ denote
its two real embeddings.  Suppose further that $p$ is inert in $F^+$.
In this setting, the framework we have established divides
$X_{U,p}(\CV)$ into four pieces:
\begin{enumerate}
\item $S=\{\tau_1,\tau_2\}$; $S^{\prime}$ is empty.  Here 
$\Sigma_{S,S^{\prime}}$ and $\Theta_{S,S^{\prime}}$ are both empty, and
thus
the scheme $Z_{S,S^{\prime}}$ is isomorphic to $X_U(\CV)$.
In particular $X_{U,p}(\CV)_{S,S^{\prime}}$ differs from $X_U(\CV)$ by
at most a Frobenius factor.  In fact, it is easy to see that in 
this case the two are isomorphic, as $X_{U,p}(\CV)_{S,S^{\prime}}$
is the image of $X_U(\CV)_{k_0}$ under the map that associates
the relative Frobenius to every abelian variety parametrized by $X_U(\CV)$.
\item $S$ is empty, $S^{\prime} = \{\tau_1,\tau_2\}$.  Again
$\Sigma_{S,S^{\prime}}$ and $\Theta_{S,S^{\prime}}$ are empty, and
the scheme $Z_{S,S^{\prime}}$ is isomorphic to $X_U(\CV)$.
Here $X_{U,p}(\CV)_{S,S^{\prime}}$ is the image of
$X_U(\CV)_{k_0}$ under the map that associates the Verschiebung
morphism to every abelian variety parametrized by $X_U(\CV)$.
\item $S = \{\tau_1\},$ $S^{\prime} = \{\tau_2\}.$  Here
$\Sigma_{S,S^{\prime}} = \{\tau_2\}$, and $\Theta_{S,S^{\prime}} = \{\tau_1\}$.
The scheme $Z^{\prime}_{S,S^{\prime}}$ is a $\PP^1 \times \PP^1$-bundle over
a zero-dimensional Shimura variety, or, more prosaically, a disjoint union
of copies of $\PP^1 \times \PP^1$'s.
\item $S = \{\tau_2\},$ $S^{\prime} = \{\tau_1\}$ exhibits the same behavior
as the previous case, with the roles of $\tau_1$ and $\tau_2$ reversed.
\end{enumerate}
It is worth noting that this description is exactly analogous to Helmuth
Stamm's description of the mod $p$ reduction of a Hilbert modular surface
with $\Gamma_0(p)$ level structure over a real quadratic field in which $p$
is inert~\cite{stamm}.  
\end{example}

\begin{example} \rm If in the previous example we had chosen to consider
a prime $p$ that was split in $F^+$, the behavior for $S = \{\tau_1,\tau_2\}$
and $S$ empty would remain the same, but when $S = \{\tau_1\}$ and 
$S^{\prime} = \{\tau_2\}$, we would have $\Theta_{S,S^{\prime}}$
and $\Sigma_{S,S^{\prime}}$ both empty.  Thus $Z_{S,S^{\prime}}$ is 
isomorphic to $X_U(\CV)$ in this case as well, and no nontrivial bundles
occur.
\end{example}

\section{A semistable model} \label{sec:semistable}

We now use Theorem~\ref{thm:DR1} to describe a semistable model for $X_U(\CV)$.
To do so it suffices to give a semistable model for $\CM$.  Note that if
$x = \CM_{\CS,\CS}$ is the ``most singular point'' of $\CM$, then the
completion $\CM_x$ is given by
$$\CM_x = W(\overline{\FF}_p)[[X_1, Y_1, \dots, X_d, Y_d]]/\<X_1Y_1 -p, \dots X_dY_d-p\>.$$
It is well-known how to resolve such a singularity, and its resolution has a
particularly nice description in terms of toric geometry.

Note first that it suffices to give a semistable model over $\ZZ[t]$
for the scheme
$$\CM^{\prime} = \spec \ZZ[t][X_1,Y_1, \dots, X_d, Y_d]/<X_1Y_1 - t, X_dY_d- t>,$$
as we can then obtain a semistable model for $\CM_x$ by base change.

The scheme $\CM^{\prime}$ can be viewed as a toric variety over $\ZZ$;
it is the $d$-fold fiber product over $\spec \ZZ[t]$ of the scheme
$$\spec \ZZ[t][X,Y]/<XY - t>.$$

As an affine toric variety (over $\ZZ$), $\spec \ZZ[t][X,Y]/<XY - t>$ corresponds
to the cone $C$ in $\RR^2$ spanned by $(1,0)$ and $(1,1)$; the morphism
to $\spec \ZZ[t]$ corresponds to the map $C \rightarrow \RR_{\geq 0}$
given by projection onto the first coordinate.  As $\CM^{\prime}$ is the $d$-fold fiber product
of this toric variety over $\spec \ZZ[t]$, it is the affine toric variety
whose cone is the $d$-fold fiber product $C_d$ of $C$ over $\RR_{\geq 0}$.  
More prosaically, $C_d$ is the cone
in $\RR^{d+1}$ over the $n$-cube $$Q = \{1\} \times [0,1] \times \dots \times [0,1].$$

Fulton (\cite{fulton}, 2.6) shows how to resolve the singularities 
of a toric variety by giving a refinement
of its fan whose cones are simplicial and unimodular.  In this case this amounts to fixing a 
decomposition $\calC$ of $Q$ into unimodular simplices, any two of which intersect
in a common face of both.  The cones over these simplices are then a fan $\CF$ consisting
of unimodular cones that refines $C_d$.  To fix ideas, we let $\calC$ be
the decomposition whose top dimensional simplices are indexed by elements of the
symmetric group $S_d$, and for a given $\sigma \in S_d$
are given by
$$C_{\sigma} = \{(1,x_1, \dots, x_d) \in Q: x_{\sigma(1)} \leq x_{\sigma(2)} \leq
\dots \leq x_{\sigma(d)}\}.$$

Let $(\CM^{\prime})^{\smallss}$ be the toric variety (over $\ZZ$) associated to $\CF$.  The map of
fans $\CF \rightarrow C_d$ induces a birational morphism from $(\CM^{\prime})^{\smallss}$ to
$\CM^{\prime}$.  As all the cones in $\CF$, the boundary of $(\CM^{\prime})^{\smallss}$ are 
unimodular, its boundary is a divisor with normal crossings.  In particular the fiber 
over $t = 0$ is a divisor with normal crossings.

After a base change from $\ZZ[t]$ to $W(\overline{\FF}_p)$, we arrive at a 
semistable resolution $\CM^{\smallss}_x$ of $\CM_x$.
The strata of the special fiber of $\CM^{\smallss}_x$ have a simple description in
terms of the combinatorics of $\calC$.  
In particular, irreducible components of
the special fiber correspond to $0$-simplices in $\calC$ (i.e., vertices of $Q$).
A collection of irreducible components intersect nontrivially if their corresponding
vertices form the vertices of a simplex in $\calC$.  For an $r$-simplex
$P$ in $\calC$, let $\CM^{\smallss}_P$ be the intersection of the irreducible components
of the special fiber of $\CM^{\smallss}_x$ that correspond to the vertices of $P$.

The image of $\CM^{\smallss}_P$ in $\CM$ is determined by the smallest face of $Q$
containing $P$.  More precisely, fix a bijection $i \mapsto \tau_i$ of $\{1, \dots, d\}$ 
with the set of real embeddings of $F^+$, and define:
$$S_P = \{\tau_i: \exists (1,x_1, \dots, x_d) \in P, x_i > 0\}$$
$$S^{\prime}_P = \{\tau_i: \exists (1,x_1, \dots, x_d) \in P, x_i < 1\}.$$
Then the image of $\CM^{\smallss}_P$ in $\CM_x$ is the intersection of
$\CM_x$ with $\CM_{S_P,S^{\prime}_P}$.  In particular, if $P$ is a top-dimensional
simplex, then $\CM^{\smallss}_P$ is a point lying over $x$.  Similarly, if
$P$ is a codimension $1$ simplex, then $\CM^{\smallss}_P$ is a rational curve; it
lies over $x$ if, and only if, $P$ is not contained in a proper face of $Q$.

The semistable resolution of $\CM_x$ gives rise to a semistable
resolution $X_U(\CV)^{\smallss}$ of $X_U(\CV)$, and the combinatorics of
$\CM^{ss}$ reflect the combinatorics of $X_U(\CV)^{\smallss}$.
In particular each vertex $V$ of $Q$ corresponds to the irreducible
component of the special fiber of $\CM^{ss}$ that dominates
$\CM_{S_V,S^{\prime}_V}$.  For each such $V$, let
$X_U(\CV)^{\smallss}_V$ denote the union of irreducible components
of the special fiber of $X_U(\CV)^{\smallss}$ that dominate some component
of $X_U(\CV)^{\smallss}$.  This is a disjoint union of smooth varieties
over $k_0$.  Extend this to simplices $P$ of arbitrary dimension in
$\calC$ by taking $X_U(\CV)^{\smallss}_P$ to be the intersection
of $X_U(\CV)^{\smallss}_V$ over all vertices $V$ of $P$.

For $0 \leq r \leq d$,
let $X_U(\CV)^{\smallss,(r)}$ denote the disjoint union of the intersections
of $r+1$ distinct irreducible components of the special fiber of
$X_U(\CV)^{\smallss}$.  As $X_U(\CV)^{\smallss}$ is semistable, these
are all smooth varieties over $\overline{\FF}_p$.

In codimensions $d$ and $d-1$ these spaces are
particularly easy to desribe.  In particular $X_U(\CV)^{\smallss,(d)}$ is simply 
the disjoint union over top-dimensional simplices $P$ in $\QQ$ of
$X_U(\CV)^{\smallss}_P$, and each of these is simply a copy of $X_U(\CV)_{\CS,\CS}$.  

Similarly, $X_U(\CV)^{\smallss,(d-1)}$ is the disjoint union of 
$X_U(\CV)^{\smallss}_P$, for $P$ a $d-1$-simplex of $\calC$.  We have:
\begin{enumerate}
\item If $P$ is contained in the relative interior of $Q$, then
$X_U(\CV)^{\smallss}_P$ is a $\PP^1$-bundle over the discrete set
$X_U(\CV)_{\CS,\CS}$. 
\item If $P$ is contained in a face $x_i = 0$ of $Q$, then
$X_U(\CV)^{\smallss}_P$ maps isomorphically onto $X_U(\CV)_{\CS \setminus \{\tau_i\},\CS}$.
\item If $P$ is contained in a face $x_i = 1$ of $Q$, then
$X_U(\CV)^{\smallss}_P$ maps isomorphically onto $X_U(\CV)_{\CS,\CS \setminus \{\tau_i\}}$.
\end{enumerate}
 
If we fix orientations on each simplex in $\calC$, we get a map
$$H^0_{\et}(X_U(\CV)^{\smallss,(d-1)}_{\overline{\FF}_p}, \QQ_l) \rightarrow
H^0_{\et}(X_U(\CV)^{\smallss,(d)}_{\overline{\FF}_p}, \QQ_l)$$
defined as follows:
For each $d$-dimensional simplex $P$, and each face $P^{\prime}$ of $P$, the map
$$H^0_{\et}(X_U(\CV)^{\smallss}_{P^{\prime}}, \QQ_l) \rightarrow H^0_{\et}(X_U(CV)^{\smallss}_P,\QQ_l)$$
is the canonical map on cohomology if the chosen orientation for $P^{\prime}$ agrees with
the orientation $P^{\prime}$ acquires from $P$, and is the negative of
this map otherwise.

We can interpret the cokernel of this map purely in terms of $X_U(\CV)$, rather
than its resolution.  In partiulcar let $X_U(\CV)^{(d-1)}$ be the disjoint
union of the $1$-strata in the special fiber of $X_U(\CV)$; that is,
the disjoint union of the strata $X_U(\CV)_{\CS,\CS \setminus \tau}$
and $X_U(\CV)_{\CS \setminus \tau, \CS}$ for all $\tau$.

We then have a commutative diagram:
$$
\begin{array}{ccc}
H^0_{\et}(X_U(\CV)^{(d-1)}_{\overline{\FF}_p}, \QQ_l)
& \rightarrow & H^0_{\et}((X_U(\CV)_{\CS,\CS})_{\overline{\FF}_p}, \QQ_l)\\
\downarrow & & \downarrow \\
H^0_{\et}(X_U(\CV)^{\smallss,(d-1)}_{\overline{\FF}_p}, \QQ_l) & \rightarrow &
H^0_{\et}(X_U(\CV)^{\smallss,(d)}_{\overline{\FF}_p}, \QQ_l) 
\end{array}
$$
where the vertical arrows are given by pullback.

An easy computation then shows:

\begin{proposition} \label{prop:combinatorial}
The vertical maps in the above diagram are a quasiisomorphism.  In
particular the cokernel of the top horizontal map is isomorphic to
the cokernel of the bottom horizontal map.
\end{proposition}

\section{Cohomology} \label{sec:cohomology}

We now illustrate how the geometric results of the previous sections, 
together with the weight spectral sequence of Rapoport-Zink,
can be used to obtain a Jacquet-Langlands correspondence on the level
of {\'e}tale cohomology.  The technical complications involved in doing
this in general are beyond the scope of this paper, and will be addressed
in future work.  Here, we will content ourselves with considering
the case where $p$ is inert in $F^+$, and $d$ is even.

Our techniques require as input Shimura varieties that are proper.
It will therefore be necessary to work not with the
moduli spaces $X_U(\CV)$, but with the more general $X_U^D(\CV)$ where
$D$ is a quaternion algebra.  We thus depart from the convention, in
force in the rest of the paper, of considering only the former.

Fix a $4$-dimensional central simple $F$-algebra $D$ split at $p$, and
a rank one left $D$-module $\CV$ equipped with
a pairing such that $r_{\tau}(\CV) = 1$ for all $\tau \in \CS$.
Consider the Shimura variety $X = X_U^D(\CV)$ for some suitable $U$.
Our first objective will be to understand the finite set 
$X_{\CS}(\overline{\FF}_p)$. 

Choose disjoint subsets $S$ and $S^{\prime}$ of $\CS$, such that
$S = \sigma S^{\prime}$.  (Since $p$ is inert in $F^+$, there is a
unique way to do this up to interchanging $S$ and $S^{\prime}$.

Fix a $\CV^{\prime}$ isomorphic to $\CV$ at all finite places 
with $r_{\tau}(\CV^{\prime}) = 2$ for $\tau \in S$
and $r_{\tau}(\CV^{\prime}) = 0$ for $\tau \in S^{\prime}$.
Thus defines a unitary group $G^{\prime}$ isomorphic to $G$ at all non-archimedean
places; let $U^{\prime}$ be the subgroup of $G^{\prime}(\AA_Q)$ corresponding to $U$.
Let $X^{\prime} = X^D_{U^{\prime}}(\CV^{\prime})$; it is a zero-dimensional
unitary Shimura variety.  We fix an identification of $\CV^{\prime}$ with
$\CV$ at all finite places.

The inclusion of $X_{\CS}$ in $X_S$, followed by
the map $X_S \rightarrow X^{\prime}$ constructed in 
section~\ref{sec:sparse}, yields a map 
$$\alpha_1: X_{\CS}(\overline{\FF}_p) \rightarrow 
X^{\prime}(\overline{\FF}_p).$$

We can construct another map from $X_{\CS}(\overline{\FF}_p)$ 
to $X^{\prime}(\overline{\FF}_p)$
as follows.  Fix a $\CV^{\prime\prime}$ isomorphic to $\CV$ at all
finite places with $r_{\tau}(\CV^{\prime\prime}) = 0$ for $\tau \in S$ and
$r_{\tau}(\CV^{\prime\prime}) = 2$ for $\tau \in S^{\prime}$.  If we let 
$X^{\prime\prime}$
denote $X_{U^{\prime\prime}}(\CV^{\prime\prime})$, and fix an
identification of $\CV^{\prime\prime}$ with $\CV$ at all finite places,
then we obtain (as in section~\ref{sec:sparse}) a map from
$X_{S^{\prime}}$ to $X^{\prime\prime}$.  On the other hand,
given $(A^{\prime\prime},\lambda^{\prime\prime},\rho^{\prime\prime}) 
\in X^{\prime\prime}(\overline{\FF}_p)$,
the $U^{\prime\prime}$-level structure $\rho^{\prime\prime}$ on 
$A^{\prime\prime}$ induces a
$U^{\prime}$-level structure on $A^{\prime\prime}$ (as we have fixed 
identifications
of $\CV$ with $\CV^{\prime}$ and $\CV^{\prime\prime}$), and
hence also a $U^{\prime}$-level structure $\rho^{\prime}$ on 
$(A^{\prime\prime})^{(p)}$.  Then 
$((A^{\prime\prime})^{(p)}, (\lambda^{\prime\prime})^{(p)}, \rho^{\prime})$ 
is a point in $X^{\prime}(\overline{\FF}_p)$.

Composing the inclusion of $X_{\CS}$ into
$X_{S^{\prime}}$ with the map $X_{S^{\prime}} \rightarrow X^{\prime\prime}$,
and then applying the above construction gives a map
$$\alpha_2: X_{\CS}(\overline{\FF}_p) \rightarrow
X^{\prime}(\overline{\FF}_p).$$

Knowing $\alpha_1(x)$ for some $x \in X_{\CS}(\overline{\FF}_p)$
is insufficient to recover $x$.  The ``missing data'' turns out to be
that of a $\Gamma_0(p)$-level structure on $\alpha_1(x)$.  Let us make
this more precise.

Let $x = (A,\lambda, \rho)$ be a point in 
$X_{\CS}(\overline{\FF}_p)$,
and let $(A^{\prime},\lambda^{\prime},\rho^{\prime})$ be the tuple
corresponding to $\alpha_1(x) \in X^{\prime}(\overline{\FF}_p)$.
The construction in section~\ref{sec:sparse} gives us an
isogeny $\pi_1:A \rightarrow A^{\prime}$ whose kernel is contained in
$A[p]$.  This
map induces a map $\CD_{A^{\prime}} \rightarrow \CD_A$, where $\CD_A$
and $\CD_{A^{\prime}}$ are the contravariant Dieudonn{\'e} modules of
$A[p]$ and $A^{\prime}[p]$ respectively.  Define: 

$$M_A \subset \bigoplus_{\tau} (\CD_A)_{p_{\tau}}$$ as follows:
\begin{enumerate}
\item $(M_A)_{p_{\tau}} = 0$ for $\tau \in S$, and
\item $(M_A)_{p_{\tau}} = V_A((\CD_A)_{p_{\sigma\tau}}) = F_A((\CD_A)_{p_{\sigma^{-1}\tau}})$
for $\tau \in S^{\prime}$.
\end{enumerate}
Note that $M_A$ is stable under $F_A$ and $V_A$.  

Let $M_{A^{\prime}}$
denote the preimage of $M_A$ under the map
$$\pi_1: \bigoplus_{\tau} (\CD_{A^{\prime}})_{p_{\tau}} \rightarrow 
\bigoplus_{\tau} (\CD_A)_{p_{\tau}}.$$
As with $M_A$, $M_{A^{\prime}}$ is stable under $F_{A^{\prime}}$ and 
$V_{A^{\prime}}$.
Note that because the kernel $L$ of $\pi_1: \CD_{A^{\prime}} \rightarrow \CD_A$
satisfies $\dim L_{p_{\tau}} = 2$ for $\tau \in S$ and $L_{p_{\tau_2}} = 0$
for $\tau \in S^{\prime}$,
we have $\dim (M_{A^{\prime}})_{p_{\tau}} = 2$
for {\em all} $\tau$ in $\CS$.
(The dimensions here are $2$ instead of $1$ because we are working with the
moduli space $X^D_U(\CV)$ rather than the simpler moduli problem $X_U(\CV)$.)

Then the subspace $M$ of $(\CD_{A^{\prime}})$ defined by
$M_{p_{\tau}} = (M_{A^{\prime}})_{p_{\tau}}$ and 
$M_{q_{\tau}} = M_{p_{\tau}}^{\perp}$ for $\tau \in \CS$ is
also stable under $F_{A^{\prime}}$ and $V_{A^{\prime}}$, and is a free 
$\OO_F/p$-module of rank $2$, stable under the action of $\OO_D/p$.
This module sits in an exact sequence
$$0 \rightarrow M \rightarrow \CD_{A^{\prime}} \rightarrow \CD(K_M) \rightarrow
0,$$
where $K_M \subset A^{\prime}[p]$ is a maximal isotropic subgroup of
$A^{\prime}[p]$ and $\CD(K_M)$ is its Dieudonn{\'e} module.

Let $B = A^{\prime}/K_M$, and let $\phi$ be the natural quotient
map.  Since $K_M$ is maximal isotropic, there is a unique prime-to-$p$
polarization $\lambda_B$ on $B$ such that $p\lambda^{\prime} =
\phi^{\vee}\lambda_B\phi$.  The cokernel of
$\phi: \CD_B \rightarrow \CD_{A^{\prime}}$ is simply $\CD_{A^{\prime}}/M$,
and is thus a locally free $\OO_F \otimes \overline{\FF}_p$-module of rank
one.
It follows that $(A^{\prime}, \lambda^{\prime}, \rho^{\prime}, B,
\lambda_B, \phi)$ is an $\overline{\FF}_p$ point on the zero-dimensional
variety $X^D_{U^{\prime},p}(\CV^{\prime})$.  
We denote this variety by $X^{\prime}_p$ in what follows.

\begin{proposition} \label{prop:level}
The map $X_{\CS}(\overline{\FF}_p) \rightarrow X^{\prime}_p(\overline{\FF}_p)$
that associates to $(A,\lambda,\rho)$ 
the tuple $(A^{\prime},\lambda^{\prime}, \rho^{\prime}, B, \lambda_B, \phi)$
constructed above is a bijection.  Moreover, under this identification
of $X_{\CS}(\overline{\FF}_p)$ with $X^{\prime}_p(\overline{\FF}_p)$,
the maps
$$\alpha_1,\alpha_2: X_{\CS} \rightarrow X^{\prime}$$
are identified with the degeneracy maps $X^{\prime}_p \rightarrow X^{\prime}$
that ``forget'' and ``mod out by'' the $\Gamma_0(p)$ level structure, respectively.
\end{proposition}
\begin{proof}
We provide an inverse to the above construction.  Let 
$(A^{\prime}, \lambda^{\prime}, \rho^{\prime},
B, \lambda_B, \phi)$ be an $\overline{\FF}_p$-point of
$X^{\prime}_p$.  We define a subspace $L$ of
$\CD_{A^{\prime}}$ by:
\begin{enumerate}
\item $L_{p_{\tau}} = \phi(\CD_B)_{p_{\tau}}$ for $\tau \in S$,
\item $L_{p_{\tau}} = 0$ for $\tau \in S^{\prime}$, and
\item $L_{q_{\tau}} = L_{p_{\tau}}^{\perp}$ for all 
$\tau \in \CS$.
\end{enumerate}
Since $F$ and $V$ are both zero map on 
$(\CD_{A^{\prime}})_{p_{\tau}}$ for $\tau \in S$,
$L$ is a Dieudonn{\'e} submodule of $\CD_{A^{\prime}}$.  Note that
if $(A^{\prime}, \lambda^{\prime}, \rho^{\prime}, B, \lambda_B, \phi)$
arose from an $(A,\lambda,\rho)$ in $X_{\{\tau_1,\tau_2\}}$ via the
above construction, then $L$ is simply the kernel of the map
$\CD_{A^{\prime}} \rightarrow \CD_A$ induced by the map $A \rightarrow
A^{\prime}$ occurring in that construction.

As usual, this module yields an exact sequence
$$0 \rightarrow L \rightarrow \CD_{A^{\prime}} \rightarrow \CD(K_L)
\rightarrow 0$$
for some maximal isotropic subgroup $K_L$ of $A^{\prime}[p]$.  Let
$A = A^{\prime}/K_L$.  The polarization $\lambda^{\prime}$ on $A^{\prime}$
induces a prime-to-$p$ polarization $\lambda$ on $A$, and the $U^{\prime}$-level
structure $\rho^{\prime}$ induces a $U$-level structure $\rho$ on $A$.
Then $(A,\lambda,\frac {1}{p}\rho)$ is a point on $X$.  Moreover,
it is clear from the previous paragraph that if 
$(A^{\prime},\lambda^{\prime},\rho^{\prime},B,\lambda_B,\phi)$
originally arose from an $\overline{\FF}_p$-point $x$ of $X_{\CS}$
via the first construction, then we have $x = (A,\lambda,\frac {1}{p}\rho)$.

We now verify that the point $(A,\lambda,\frac {1}{p}\rho)$ we have
constructed lies on $X_{\CS}$.
Consider the Dieudonn{\'e} modules $\tCD_A$, $\tCD_B$, and
$\tCD_{A^{\prime}}$ of $A[p^{\infty}]$, $B[p^{\infty}]$, and 
$A^{\prime}[p^{\infty}]$ respectively.  Let $\varphi$ denote the
quotient map $A^{\prime} \rightarrow A$.

The map $\varphi$ identifies $\tCD_A$ with a submodule
of $\tCD_{A^{\prime}}$, and in particular
identifies $(\tCD_A)_{p_{\tau}}$ with $p(\tCD_{A^{\prime}})_{p_{\tau}}$
for $\tau \in S^{\prime}$.
For $\tau \in S$, the maps $F_{A^{\prime}}$ and $V_{A^{\prime}}$  
map $(\tCD_{A^{\prime}})_{p_{\sigma^{-1}\tau}}$ 
and $(\tCD_{A^{\prime}})_{p_{\sigma\tau}}$ isomorphically
onto $(\tCD_{A^{\prime}})_{p_{\tau}}$, and so map $(\tCD_A)_{p_{\sigma^{-1}\tau}}$
and $(\tCD_A)_{P_{\sigma\tau}}$
isomorphically onto $p(\tCD_{A^{\prime}})_{p_{\tau}}$.  In particular
$F_{A^{\prime}}(\tCD_A)_{p_{\sigma^{-1}\tau}} = V_{A^{\prime}}(\tCD_A)_{p_{\sigma\tau}}$ 
and so 
$(A,\lambda,\frac {1}{p}\rho)$ lies in $X_S$.

On the other hand, by our construction of $A$ we have
$\phi(\tCD_B)_{p_{\tau}} = \varphi(\tCD_A)_{p_{\tau}}$ for $\tau \in S$.  But
$F_B$ and $V_B$ map $(\tCD_B)_{p_{\sigma^{-1}\tau}}$ and
$(\tCD_B)_{p_{\sigma\tau}}$ (respectively) isomorphically onto
$(\tCD_B)_{p_{\tau}}$.  

We therefore have:
$$F_{A^{\prime}}\phi(\tCD_B)_{p_{\tau}} = 
F_{A^{\prime}}\phi(V_B(\tCD_B)_{p_{\sigma\tau}}
= F_{A^{\prime}}V_{A^{\prime}}\phi(\tCD_B)_{p_{\sigma\tau}} =
p\phi(\tCD_B)_{p_{\sigma\tau}}.$$
Since $\sigma^2\tau$ is also in $S$, we similarly have
$$V_{A^{\prime}}\phi(\tCD_B)_{p_{\sigma^2\tau}} =
p\phi(\tCD_B)_{p_{\sigma\tau}}.$$
It follows that $F_A(\tCD_A)_{p_{\tau}} = 
V_A(\tCD_A)_{p_{\sigma^2\tau}},$ for all $\tau \in S$,
so $(A,\lambda, \frac {1}{p}\rho)$ lies in $X_{S^{\prime}}$.

One now verifies easily that the construction given above is both
a left and right inverse to the map $X_{\CS}(\overline{\FF}_p)
\rightarrow X^{\prime}_p(\overline{\FF}_p)$. 

A straightforward calculation shows that under this identification, the maps
$\alpha_1, \alpha_2$ from $X_{\CS}$ to
$X^{\prime}$ correspond exactly to the two degeneracy maps
$X^{\prime}_p \rightarrow X^{\prime}$ that
``forget the level structure at $p$'' and ``mod out by the level structure
at $p$'' respectively.
\end{proof}

\begin{remark} \rm Although we will make no use of the fact here, it is not
difficult to generalize this construction to a description of the
dimension zero locus for any totally real field $F^+$ in which every
prime of $F^+$ over $p$ has {\em even} residue class degree.  When primes
of odd residue class degree occur over $p$, the construction here breaks
down completely.  It is likely that in this case one needs to consider
an auxiliary quaternion algebra {\em ramified} at the primes of odd
residue class degree over $p$.  Even if one does this, it is not clear
at all what archimedean invariants $r_{\tau}^{\prime}, s_{\tau}^{\prime}$
should be taken.
\end{remark} 

We now turn to arithmetic questions.  The key tool here will be the
weight spectral sequence of Rapoport-Zink~\cite{RZ2},~\cite{Saito}.  
This applies
only to proper schemes with semistable reduction.  In order to obtain
a proper scheme to apply it on, we will need to make the assumption
that $D$ is ramified at at least one prime- that is, $D$ is a quaternion
algebra rather than a matrix algebra.

In this case, the scheme $X^D_{U,p}(\CV)$ is proper over $W(k_0)$, 
and has the semistable resolution $X^{\smallss}$ described in 
section~\ref{sec:semistable}.  For each $r$, with $0 \leq r \leq d$,
let $Y^{(r)}$ be the scheme over $\overline{\FF}_p$ consisting
of the disjoint union of all nonempty $r+1$-fold intersections
of irreducible components of the special fiber of $X^{\smallss}$.

Fix a prime $l$ different from $p$.  The weight spectral sequence 
(c.f.~\cite{Saito}, Corollary 2.2.4) is a spectral sequence
$$E_1^{p,q} = \bigoplus_{i \geq \max(0,-p)} H^{q-2i}_{\et}(Y^{(p+2i)}, \QQ_l(-i))
\Rightarrow H^{p+q}_{\et}(X^{\smallss}_{\overline{K}}, \QQ_l).$$
This spectral sequence degenerates at $E_2$.  In our situation, its nonzero
terms are concentrated between $p= -d$ and $p=d$.  Thus the sequence
$$E_2^{d,0}, E_2^{d-1,1}, \dots E_2^{1-d,2d-1}, E_2^{-d,2d}$$
is the sequence of successive quotients of a $2d+1$-step filtration
(the ``weight filtration'') on $H^d_{\et}(X^{\smallss}_{\overline{K}}, \QQ_l)$.
(Equivalently, this is a filtration of
$H^d_{\et}(X^D_{U,p}(\CV)_{\overline{K}})$, as the blow-up we
performed altered only the special fiber.)

We will compute the top and bottom quotients of this filtration.  
First we need some notation:
Let $U^{\prime}_p$ be the subgroup of $G^{\prime}(\AA_{\QQ})$ equal to
$U^{\prime}$ at all prime to $p$ places but isomorphic to $\Gamma_0(p)$
at $p$.
Let $S_2(G^{\prime},U^{\prime}_p,\QQ_l)$ denote the space of algebraic modular 
forms~\cite{gross} on $G^{\prime}$ of level $U^{\prime}_p$ and ``weight 2''; that is, the space of
$\QQ_l$-valued functions on the double coset space $G^{\prime}(\QQ)\backslash G^{\prime}(\AA_{\QQ})/U^{\prime}_p$.
There are two natural degeneracy maps
$$S_2(G^{\prime},U^{\prime},\QQ_l) \rightarrow S_2(G^{\prime},U^{\prime}_p,\QQ_l);$$
let $S_2(G^{\prime},U^{\prime}_p,\QQ_l)^{\pnew}$ denote the ``$p$-new'' quotient of
$S_2(G^{\prime},U^{\prime}_p,\QQ_l)$ by the images of the two degeneracy maps.
Since spaces of functions on a finite set are naturally self-dual, we also
have a $p$-new subspace $S_2(G^{\prime},U^{\prime}_p,\QQ_l)_{\pnew}$ defined as the
intersection of the kernels of the duals of the degeneracy maps.

Note that as finite sets, $X^{\prime}_p(\overline{\FF}_p)$ and $X^{\prime}(\overline{\FF}_p)$
are in bijection with the double coset spaces $G^{\prime}(\QQ)\backslash G^{\prime}(\AA_{\QQ})/U^{\prime}_p$
and $G^{\prime}(\QQ)\backslash G^{\prime}(\AA_{\QQ})/U^{\prime}$, respectively.  These bijections
identify the degeneracy maps on $S_2(G^{\prime},U^{\prime},\QQ_l)$ with the pullback
maps
$$H^0_{\et}((X^{\prime})_{\overline{\FF}_p}, \QQ_l) 
\rightarrow H^0_{\et}((X^{\prime}_p)_{\overline{\FF}_p}, \QQ_l)$$
induced by the two degeneracy maps $X^{\prime}_p \rightarrow X^{\prime}$.

\begin{theorem} \label{thm:filtration}
The highest and lowest weight quotients of
the weight filtration on 
$H^d((X^D_{U,p})_{\overline{K}}, \QQ_l)$ are naturally isomorphic to 
$S_2(G^{\prime},U^{\prime}_p,\QQ_l)_{\pnew}$ and
$S_2(G^{\prime},U^{\prime}_p,\QQ_l)^{\pnew},$ respectively.

Moreover, these isomorphisms are compatible with the actions of the prime-to-$p$
Hecke operators for $G$ and $G^{\prime}$, respectively.
\end{theorem}
\begin{proof}
We prove this statement for the lowest weight quotient; the corresponding statement for the
top then follows by the fact that the weight spectral sequence is self-dual under
Poincar{\'e} duality.

The bottom quotient of the weight filtration is the term $E_2^{d,0}$ of the weight spectral
sequence; that is, the cokernel of the map $E_1^{d-1,0} \rightarrow E_1^{d,0}$.
This is the restriction map: 
$$H^0_{\et}(Y^{(d-1)},\QQ_l) \rightarrow H^0_{\et}(Y^{(d)},\QQ_l).$$

By Proposition~\ref{prop:combinatorial} the cokernel of this restriction map is
naturally isomorphic to the cokernel of the map
$$\bigoplus_{\tau} H^0_{\et}((X_{\CS \setminus \{\tau\},\CS})_{\overline{\FF}_p}, \QQ_l)
\oplus H^0_{\et}((X_{\CS,\CS \setminus \{\tau\}})_{\overline{\FF}_p}, \QQ_l) \rightarrow
H^0_{\et}((X_{\CS,\CS})_{\overline{\FF}_p}, \QQ_l).$$
For any given $\tau$, either $X_{\CS \setminus \{\tau\}, \CS}$ is contained in $X_{S,S^{\prime}}$
and $X_{\CS, \CS \setminus \{\tau\}}$ is contained in $X_{S^{\prime},S}$, or vice versa.

Suppose $\tau \in S^{\prime}$, so $X_{\CS \setminus \{\tau\}, \CS}$ is contained in $X_{S,S^{\prime}}$.
The scheme $X_{S,S^{\prime}}$ is a disjoint union of copies of $(\PP^1)^d$, indexed by the
points of $X^{\prime}$; the intersection
of $X_{\CS \setminus \tau, \CS}$ with each of these copies of $(\PP^1)^d$ is a rational curve.
In particular pullback induces an isomorphism:
$$H^0_{\et}((X_{\CS \setminus \tau, \CS})_{\overline{\FF}_p}, \QQ_l) \cong
H^0_{\et}((X_{S,S^{\prime}})_{\overline{\FF}_p}, \QQ_l) \cong 
H^0_{\et}(X^{\prime}_{\overline{\FF}_p}, \QQ_l).$$
If we make these identifications, and we identify $X_{\CS,\CS}(\overline{\FF}_p)$ with $X^{\prime}_p(\overline{\FF}_p)$
via Proposition~\ref{prop:level}, the pullback map 
$$H^0_{\et}((X_{\CS \setminus \{\tau\}, \CS})_{\overline{\FF}_p}, \QQ_l) \rightarrow
H^0_{\et}((X_{\CS,\CS})_{\overline{\FF}_p}, \QQ_l)$$
becomes the degeneracy map
$$H^0_{\et}((X^{\prime})_{\overline{\FF}_p}, \QQ_l) \rightarrow
H^0_{\et}((X^{\prime}_p)_{\overline{\FF}_p}, \QQ_l)$$
that comes from ``forgetting the level structure at p.''

Similarly, the pullback map
$$H^0_{\et}((X_{\CS, \CS \setminus \{\tau\}})_{\overline{\FF}_p}, \QQ_l) \rightarrow
H^0_{\et}((X_{\CS,\CS})_{\overline{\FF}_p}, \QQ_l)$$
becomes the degeneracy map
$$H^0_{\et}((X^{\prime})_{\overline{\FF}_p}, \QQ_l) \rightarrow
H^0_{\et}((X^{\prime}_p)_{\overline{\FF}_p}, \QQ_l)$$
that ``mods out by the level structure at p.''

If $\tau$ lies in $S$ rather than $S^{\prime}$, we still get the above two degeneracy
maps, but in the opposite order.

Thus $E_2^{d,0}$ is isomorphic to the quotient of
$H^0_{\et}((X^{\prime}_p)_{\overline{\FF}_p}, \QQ_l)$ by the image
of the two degeneracy maps from $H^0_{\et}(X^{\prime}_{\overline{\FF}_p}, \QQ_l).$
By the above discussion, this is isomorphic to 
$S_2(G^{\prime}, U^{\prime}_p, \QQ_l)^{\pnew}.$

For the Hecke equivariance, note that we have
fixed an isomorphism between $G$ and $G^{\prime}$ at all finite places,
and therefore can associate to every prime-to-$p$ Hecke operator
for $G$ a corresponding prime-to-$p$ Hecke operator for $G^{\prime}$.
Such Hecke operators act via correspondences on $X_{U,p}$
or $X^{\prime}_p$ defined in terms of the prime-to-$p$ level structures
of the points of these varieties.

Since every construction in this paper is defined in terms of $p$-isogenies,
these constructions commute with correspondences defined in this way, and
are thus compatible with the Hecke actions in the above sense.
\end{proof}

\begin{remark} \rm The above construction can be easily adapted
to work with coefficients in an arbitrary lisse sheaf over 
$\overline{\QQ}_l$.  Thus the result above extends to ``higher weight''.
In certain circumstances one can even obtain results with
coefficients in $\ZZ_l$; the main difficulty is that
the weight spectral sequence is not known to degenerate at $E_2$ in this
case.  This can be checked by hand in many cases, such as when $d=2$.
\end{remark}

A special case of Jacquet-Langlands for $G$ and $G^{\prime}$
follows:

\begin{corollary} Under the assumptions of this section,
let $\pi$ be an automorphic representation of $G$ (resp. $G^{\prime}$) that is 
cohomological, and Steinberg at a prime $p$ split in $E$
and inert in $F^+$.  Then there is an automorphic representation
$\pi^{\prime}$ of $G^{\prime}$ (resp. $G$) isomorphic to $\pi$ at
all non-archimedean places.
\end{corollary}

\section{Appendix: Bilinear pairings} \label{sec:appendix}

Let $E$, $F^+$, and $F$ be as in section~\ref{sec:basic}, and let
$D$ be a central simple $F$-algebra of dimension $n^2$, equipped with
an involution $\alpha \mapsto \overline{\alpha}$ of the second kind.
Let $\CV$ be a free left $D$-module of rank one, equipped with a pairing
$\langle,\rangle: V \times V \rightarrow \QQ$ that is alternating and
nondegenerate, and satisfies $\langle \alpha x, y \rangle = 
\langle x, \overline{\alpha} y \rangle$.

Note that $D$ is split at infinity, as $F$ is purely imaginary.
For each real embedding $\tau$ of $F^+$, $D \otimes_{F^+,\tau} \RR$
is isomorphic to $M_n(\CC)$.  For a suitable choice of idempotent
$e$ in $D \otimes_{F^+,\tau} \RR$ such that $e = \overline{e}$,
$e(\CV \otimes_{F^+,\tau} \RR)$ is an $n$-dimensional complex vector
space, on which $\langle,\rangle$ induces a nondegenerate Hermitian
pairing.  We let $r_{\tau}(\CV)$ and $s_{\tau}(\CV)$ denote the number
of $1$'s and $-1$'s in the signature of this pairing; this is independent
of the particular idempotent chosen and satisfies 
$r_{\tau}(\CV) + s_{\tau}(\CV) = n$. 

We are concerned with the question of determining whether, for given
$r_{\tau}^{\prime}$ and $s_{\tau}^{\prime}$, there exists a
$\CV^{\prime}$ and $\langle,\rangle^{\prime}$ such that the pair
$(\CV^{\prime}, \langle,\rangle^{\prime})$ is isomorphic to
$(\CV, \langle, \rangle)$ over $\AA^f_{\QQ}$ but such that
$r_{\tau}(\CV^{\prime}) = r_{\tau}^{\prime}$ and $s_{\tau}(\CV^{\prime})
= s_{\tau}^{\prime}$.  This turns out to be an easy application of
results of Kottwitz~\cite{Ko}.

Let $G$ be the subgroup of $\Aut_D(\CV)$ of automorphisms preserving
$\langle,\rangle$.  (Note that this is different from the $G$ appearing
in the main body of the text, as we are not working up to a scalar multiple.)
The cohomology set $H^1(\QQ,G)$ is in bijection with the set of
left $D$-modules $\CV^{\prime}$ of rank one with pairings 
$\langle,\rangle^{\prime}$ of the sort under consideration.  A pair
$(\CV^{\prime},\langle,\rangle^{\prime})$ is isomorphic to
$(\CV,\langle,\rangle)$ over $\AA^f_{\QQ}$ if and only if the corresponding 
class in $H^1(Q_v,G)$ is trivial for all finite places $v$ of $\QQ$.

It thus suffices to understand the image of the natural map
$$H^1(Q,G) \rightarrow \oplus_v H^1(Q_v,G)$$
where $v$ runs over the places of $\QQ$.  Let $A(G)$ 
denote $\pi_0(Z(\hat G)^{\gal(\overline{\QQ}/\QQ)})^{\vee}$,
i.e., the Pontryagin dual of group of components of the Galois invariants
of the center of the connected Langlands dual of $G$.  Similarly let
$A(G_v)$ denote $\pi_0(Z(\hat G_v)^{\gal(\overline{\QQ_v}/\QQ_v)})^{\vee}$.
Then Kottwitz~\cite{Ko} constructs for each $v$ a natural map
$$H^1(Q_v,G) \rightarrow A(G_v),$$
and shows (\cite{Ko}, Proposition 2.6) that the image of
$H^1(Q,G)$ in $\oplus_v H^1(Q_v,G)$ is equal to the kernel of the
composite map
$$\oplus_v H^1(Q_v, G) \rightarrow \oplus_v A(G_v) \rightarrow A(G).$$

In light of this, understanding our question amounts to understanding
the map $A(G_{\infty}) \rightarrow A(G)$.  The group $G \times_{\QQ} \RR$
is isomorphic to a product
$$\prod_{\tau \in \CS} U(r_{\tau},s_{\tau})$$ of unitary
groups.

Over $\CC$, a unitary group $U(r,s)$ becomes isomorphic to $\GL_n(\CC)$,
and comes equipped with a natural action of $\gal(\CC/\RR)$.
Complex conjugation acts on this $\GL_n(\CC)$ by sending a matrix
$g$ to $(g^*)^{-1}$, where $*$ denotes adjoint with respect to a Hermitian
pairing on $\CC^n$ with signature consisting of $r$ $1$'s and $s$ $-1$'s.
If we fix a maximal torus of $T$ of $U(r,s)$, the induced action on
the character group of $T$ sends a character to its inverse.

This action of $\gal(\CC/\RR)$ on the character group of $T$ induces
an action of $\gal(\CC/\RR)$ on the root datum of $\GL_n(\CC)$, and
hence also on the dual root datum.  This action in turn defines a unique
action of $\gal(\CC/\RR)$ on the Langlands dual of $\GL_n(\CC)$ that
respects the complex structure of $\GL_n(\CC)$.  (Note that by contrast
the Galois action on $U(r,s) \times_{\RR} \CC$ did not respect this
complex structure.)  
Thus the connected Langlands dual of $U(r,s)$ is again $\GL_n(\CC)$,
equipped with an action of $\gal(\CC/\RR)$.  The action of complex
conjugation on the maximal torus of this group can be read off of the action
on the root datum; it sends any element of the maximal torus to its inverse.
In particular, $Z(\widehat U(r,s))^{\gal(\CC/\RR)}$ is equal to
$$\{x \in \CC: x = x^{-1}\} = \{\pm 1\}.$$  Moreover, under the
correspondence between $H^1(\RR, U(r,s))$ and isomorphism classes
of nondegenerate Hermitian forms on $\CC^{r+s}$, the element of
$H^1(\RR, U(r,s))$ corresponding to a Hermitian form of signature
consisting of $r^{\prime}$ $1$'s and $s^{\prime}$ $-1$'s maps
to the element $(-1)^{r-r^{\prime}}$ under the map
$H^1(\RR, U(r,s)) \rightarrow A(U(r,s))$.

Thus we have a natural isomorphism:
$$A(G_{\infty}) \rightarrow \prod_{\tau \in \CS} \{\pm 1\}.$$
On the other hand, we can compute $A(G)$ in exactly the same fashion.
In particular $G \times_{\QQ} E$ is isomorphic to $D^{\times}$, and
$\gal(E/\QQ)$ acts on $B^{\times}$ by sending $\alpha$ to
$\overline{\alpha}^{-1}$.  Thus the center of $G \times_{\QQ} \overline{\QQ}$
is isomorphic to $(F \otimes_E \overline{\QQ})^{\times}$, where 
$\gal(\overline{\QQ}/\QQ)$ acts on $F$ via the action of $\gal(E/\QQ)$
on $B$ and on $\overline{\QQ}$ in the obvious way.

If we fix a maximal torus $T$ of $G$, we obtain
a sequence of maps
$$\chi(Z(\hat G))^* \rightarrow \chi(\hat T)^* = \chi(T) \rightarrow
\chi(Z(G))$$
where $\chi$ denotes the character group.  This identifies $\chi(Z(\hat G))^*$
with a sublattice of $\chi(Z(G))$ in a manner that respects the action
of $\gal(\overline{\QQ}/\QQ)$.  Moreover, because over $\overline{\QQ}$,
$G$ is simply a product of copies of $\GL_n$, the image of this
sequence of maps is just $n\chi(Z(G))$.  It follows that $Z(\hat G)$
is simply $(\CC^{\times})^{\CS}$, on which
$\gal(F^+/\QQ)$ acts by permuting the factors and complex conjugation
in $\gal(E/\QQ)$ acts by inversion.  In particular 
$Z(\hat G)^{\gal(\overline{\QQ}/\QQ)}$ is simply ${\pm 1}$.  The map
$A(G_{\infty}) \rightarrow A(G)$ simply takes an element of
$\prod_{\tau \in \CS} \{\pm 1\}$ and multiplies the entries
together.

The upshot of this is:
\begin{proposition} \label{prop:pairing}
Let $r_{\tau}^{\prime}$ and $s_{\tau}^{\prime}$ be integers between
$0$ and $n$ such that $s_{\tau}^{\prime} = n - r_{\tau}^{\prime}$ for
all $\tau$.  
Then there exists a $(\CV^{\prime},\langle,\rangle^{\prime})$
isomorphic to $(\CV,\langle,\rangle)$ at all finite places but with
archimedean invariants $r_{\tau}^{\prime},s_{\tau}^{\prime}$ if and only
if
$$\sum_{\tau} r_{\tau} \equiv \sum_{\tau} r_{\tau}^{\prime} \text{mod $2$}.$$
\end{proposition}
\begin{proof}
Let $x$ be the class in $\oplus_v H^1(\QQ_v,G)$ that is trivial at all
finite places but that corresponds to $r_{\tau}^{\prime},s_{\tau}^{\prime}$
at infinity.  Then $x$ maps to 
$(-1)^{\sum_{\tau} r_{\tau} - r_{\tau}^{\prime}}$ in $A(G)$, and is thus
in the kernel of the map to $A(G)$ if and only if the parity condition
given above is satisifed.  Since the kernel of this map is equal to
the image of $H^1(Q,G)$, the result follows immediately.
\end{proof}

\begin{corollary} \label{cor:pairing}
Let $\CV$ be an $n$-dimensional $F$-vector space, together with an
althernating nondegenerate pairing $\langle,\rangle: \CV \times \CV \rightarrow
\QQ$, such that $\langle \alpha x, y\rangle = \langle x, \overline{\alpha} y
\rangle$.  Let $r_{\tau}^{\prime}$ and $s_{\tau}^{\prime}$ be
integers between $0$ and $n$ such that $s_{\tau}^{\prime} = 
n - r_{\tau}^{\prime}$ for all $\tau$.
Then there exists a $(\CV^{\prime},\langle,\rangle^{\prime})$
isomorhpic to $(\CV,\langle,\rangle)$ at all finite places but with
archimedean invariants $r_{\tau}^{\prime},s_{\tau}^{\prime}$ if and only
if
$$\sum_{\tau} r_{\tau} \equiv \sum_{\tau} r_{\tau}^{\prime} \text{mod $2$}.$$
\end{corollary}
\begin{proof}
Let $D = \GL_n(F)$, and let $\alpha \mapsto \overline{\alpha}$ be
the conjugate transpose on $D$.  Then to any $(\CV,\langle,\rangle)$
as in the statement of the corollary, we can associate the left $D$-module
$\CV^n$, together with the natural pairing induced by $\langle,\rangle$.  
This is an equivalence of categories between $(\CV,\langle,\rangle)$ as in
the statement of the corollary and those considered in the preceding
proposition (for $D = GL_n(F)$).  Thus the corollary follows immediately
from the preceding proposition.
\end{proof}

\textsc{Acknowledgements.}
The author is grateful to Richard Taylor for his continued interest and
helpful suggestions, to Johan de Jong for his help with the deformation
theory of abelian varieties, to Brian Osserman for many
helpful conversations, and to Marc Hubert-Nicole for his helpful comments on
an earlier version of this paper.  The author was partially supported by 
the National Science Foundation. 

\providecommand{\bysame}{\leavevmode\hbox to3em{\hrulefill}\thinspace}

\end{document}